\documentstyle[fullpage,12pt]{article}
\date{June 1996}
\input tcilatex

\makeatletter

\def\diagram{\leftwidth=\z@ \rightwidth=\z@ \topheight=\z@
\botheight=\z@ \setbox\@picbox\hbox\bgroup}

\def\enddiagram{\egroup\wd\@picbox\rightwidth\unitlength
\ht\@picbox\topheight\unitlength \dp\@picbox\botheight\unitlength
\hskip\leftwidth\unitlength\box\@picbox}

\def\bfig{\begin{diagram}}
\def\efig{\end{diagram}}
\newcount\wideness \newcount\leftwidth \newcount\rightwidth
\newcount\highness \newcount\topheight \newcount\botheight

\def\ratchet#1#2{\ifnum#1<#2 \global #1=#2 \fi}

\def\putbox(#1,#2)#3{%
\horsize{\wideness}{#3} \divide\wideness by 2
{\advance\wideness by #1 \ratchet{\rightwidth}{\wideness}}
{\advance\wideness by -#1 \ratchet{\leftwidth}{\wideness}}
\vertsize{\highness}{#3} \divide\highness by 2
{\advance\highness by #2 \ratchet{\topheight}{\highness}}
{\advance\highness by -#2 \ratchet{\botheight}{\highness}}
\put(#1,#2){\makebox(0,0){$#3$}}}

\def\putlbox(#1,#2)#3{%
\horsize{\wideness}{#3}
{\advance\wideness by #1 \ratchet{\rightwidth}{\wideness}}
{\ratchet{\leftwidth}{-#1}}
\vertsize{\highness}{#3} \divide\highness by 2
{\advance\highness by #2 \ratchet{\topheight}{\highness}}
{\advance\highness by -#2 \ratchet{\botheight}{\highness}}
\put(#1,#2){\makebox(0,0)[l]{$#3$}}}

\def\putrbox(#1,#2)#3{%
\horsize{\wideness}{#3}
{\ratchet{\rightwidth}{#1}}
{\advance\wideness by -#1 \ratchet{\leftwidth}{\wideness}}
\vertsize{\highness}{#3} \divide\highness by 2
{\advance\highness by #2 \ratchet{\topheight}{\highness}}
{\advance\highness by -#2 \ratchet{\botheight}{\highness}}
\put(#1,#2){\makebox(0,0)[r]{$#3$}}}

\def\adjust[#1]{} 

\newcount \coefa
\newcount \coefb
\newcount \coefc
\newcount\tempcounta
\newcount\tempcountb
\newcount\tempcountc
\newcount\tempcountd
\newcount\xext
\newcount\yext
\newcount\xoff
\newcount\yoff
\newcount\gap%
\newcount\arrowtypea
\newcount\arrowtypeb
\newcount\arrowtypec
\newcount\arrowtyped
\newcount\arrowtypee
\newcount\height
\newcount\width
\newcount\xpos
\newcount\ypos
\newcount\run
\newcount\rise
\newcount\arrowlength
\newcount\halflength
\newcount\arrowtype
\newdimen\tempdimen
\newdimen\xlen
\newdimen\ylen
\newsavebox{\tempboxa}%
\newsavebox{\tempboxb}%
\newsavebox{\tempboxc}%

\newdimen\w@dth

\def\setw@dth#1#2{\setbox\z@\hbox{$#1$}\w@dth=\wd\z@
\setbox\@ne\hbox{$#2$}\ifnum\w@dth<\wd\@ne \w@dth=\wd\@ne \fi
\advance\w@dth by 1.2em}

\def\t@^#1_#2{\def\n@one{#1}\def\n@two{#2}\mathrel{\setw@dth{#1}{#2}
\mathop{\hbox to \w@dth{\rightarrowfill}}\limits
\ifx\n@one\empty\else ^{\box\z@}\fi
\ifx\n@two\empty\else _{\box\@ne}\fi}}
\def\t@@^#1{\@ifnextchar_ {\t@^{#1}}{\t@^{#1}_{}}}
\def\to{\@ifnextchar^ {\t@@}{\t@@^{}}}

\def\t@left^#1_#2{\def\n@one{#1}\def\n@two{#2}\mathrel{\setw@dth{#1}{#2}
\mathop{\hbox to \w@dth{\leftarrowfill}}\limits
\ifx\n@one\empty\else ^{\box\z@}\fi
\ifx\n@two\empty\else _{\box\@ne}\fi}}
\def\t@@left^#1{\@ifnextchar_ {\t@left^{#1}}{\t@left^{#1}_{}}}
\def\toleft{\@ifnextchar^ {\t@@left}{\t@@left^{}}}

\def\two@^#1_#2{\def\n@one{#1}\def\n@two{#2}\mathrel{\setw@dth{#1}{#2}
\mathop{\vcenter{\hbox to \w@dth{\rightarrowfill}\kern-1.7ex
                 \hbox to \w@dth{\rightarrowfill}}%
       }\limits
\ifx\n@one\empty\else ^{\box\z@}\fi
\ifx\n@two\empty\else _{\box\@ne}\fi}}
\def\tw@@^#1{\@ifnextchar_ {\two@^{#1}}{\two@^{#1}_{}}}
\def\two{\@ifnextchar^ {\tw@@}{\tw@@^{}}}

\def\tofr@^#1_#2{\def\n@one{#1}\def\n@two{#2}\mathrel{\setw@dth{#1}{#2}
\mathop{\vcenter{\hbox to \w@dth{\rightarrowfill}\kern-1.7ex
                 \hbox to \w@dth{\leftarrowfill}}%
       }\limits
\ifx\n@one\empty\else ^{\box\z@}\fi
\ifx\n@two\empty\else _{\box\@ne}\fi}}
\def\t@fr@^#1{\@ifnextchar_ {\tofr@^{#1}}{\tofr@^{#1}_{}}}
\def\tofro{\@ifnextchar^ {\t@fr@}{\t@fr@^{}}}

\def\mon{\mathop{\m@th\hbox to
      14.6\P@{\lasyb\char'51\hskip-2.1\P@$\arrext$\hss
$\mathord\rightarrow$}}\limits} 
\def\leftmono{\mathrel{\m@th\hbox to
14.6\P@{$\mathord\leftarrow$\hss$\arrext$\hskip-2.1\P@\lasyb\char'50%
}}\limits} 
\mathchardef\arrext="0200       

\def\settypes(#1,#2,#3){\arrowtypea#1 \arrowtypeb#2 \arrowtypec#3}
\def\settoheight#1#2{\setbox\@tempboxa\hbox{#2}#1\ht\@tempboxa\relax}%
\def\settodepth#1#2{\setbox\@tempboxa\hbox{#2}#1\dp\@tempboxa\relax}%
\def\settokens[#1`#2`#3`#4]{%
     \def\tokena{#1}\def\tokenb{#2}\def\tokenc{#3}\def\tokend{#4}}
\def\setsqparms[#1`#2`#3`#4;#5`#6]{%
\arrowtypea #1
\arrowtypeb #2
\arrowtypec #3
\arrowtyped #4
\width #5
\height #6
}
\def\setpos(#1,#2){\xpos=#1 \ypos#2}

\def\settriparms[#1`#2`#3;#4]{\settripairparms[#1`#2`#3`1`1;#4]}%

\def\settripairparms[#1`#2`#3`#4`#5;#6]{%
\arrowtypea #1
\arrowtypeb #2
\arrowtypec #3
\arrowtyped #4
\arrowtypee #5
\width #6
\height #6
}

\def\resetparms{\settripairparms[1`1`1`1`1;500]\width 500}

\resetparms

\def\mvector(#1,#2)#3{
\put(0,0){\vector(#1,#2){#3}}%
\put(0,0){\vector(#1,#2){26}}%
}
\def\evector(#1,#2)#3{{
\arrowlength #3
\put(0,0){\vector(#1,#2){\arrowlength}}%
\advance \arrowlength by-30
\put(0,0){\vector(#1,#2){\arrowlength}}%
}}

\def\horsize#1#2{%
\settowidth{\tempdimen}{$#2$}%
#1=\tempdimen
\divide #1 by\unitlength
}

\def\vertsize#1#2{%
\settoheight{\tempdimen}{$#2$}%
#1=\tempdimen
\settodepth{\tempdimen}{$#2$}%
\advance #1 by\tempdimen
\divide #1 by\unitlength
}

\def\putvector(#1,#2)(#3,#4)#5#6{{%
\ifnum3<\arrowtype
\putdashvector(#1,#2)(#3,#4)#5\arrowtype
\else
\ifnum\arrowtype<-3
\putdashvector(#1,#2)(#3,#4)#5\arrowtype
\else
\xpos=#1
\ypos=#2
\run=#3
\rise=#4
\arrowlength=#5
\ifnum \arrowtype<0
    \ifnum \run=0
        \advance \ypos by-\arrowlength
    \else
        \tempcounta \arrowlength
        \multiply \tempcounta by\rise
        \divide \tempcounta by\run
        \ifnum\run>0
            \advance \xpos by\arrowlength
            \advance \ypos by\tempcounta
        \else
            \advance \xpos by-\arrowlength
            \advance \ypos by-\tempcounta
        \fi
    \fi
    \multiply \arrowtype by-1
    \multiply \rise by-1
    \multiply \run by-1
\fi
\ifcase \arrowtype
\or \put(\xpos,\ypos){\vector(\run,\rise){\arrowlength}}%
\or \put(\xpos,\ypos){\mvector(\run,\rise)\arrowlength}%
\or \put(\xpos,\ypos){\evector(\run,\rise){\arrowlength}}%
\fi\fi\fi
}}

\def\putsplitvector(#1,#2)#3#4{
\xpos #1
\ypos #2
\arrowtype #4
\halflength #3
\arrowlength #3
\gap 140
\advance \halflength by-\gap
\divide \halflength by2
\ifnum\arrowtype>0
   \ifcase \arrowtype
   \or \put(\xpos,\ypos){\line(0,-1){\halflength}}%
       \advance\ypos by-\halflength
       \advance\ypos by-\gap
       \put(\xpos,\ypos){\vector(0,-1){\halflength}}%
   \or \put(\xpos,\ypos){\line(0,-1)\halflength}%
       \put(\xpos,\ypos){\vector(0,-1)3}%
       \advance\ypos by-\halflength
       \advance\ypos by-\gap
       \put(\xpos,\ypos){\vector(0,-1){\halflength}}%
   \or \put(\xpos,\ypos){\line(0,-1)\halflength}%
       \advance\ypos by-\halflength
       \advance\ypos by-\gap
       \put(\xpos,\ypos){\evector(0,-1){\halflength}}%
   \fi
\else \arrowtype=-\arrowtype
   \ifcase\arrowtype
   \or \advance \ypos by-\arrowlength
       \put(\xpos,\ypos){\line(0,1){\halflength}}%
       \advance\ypos by\halflength
       \advance\ypos by\gap
       \put(\xpos,\ypos){\vector(0,1){\halflength}}%
   \or \advance \ypos by-\arrowlength
       \put(\xpos,\ypos){\line(0,1)\halflength}%
       \put(\xpos,\ypos){\vector(0,1)3}%
       \advance\ypos by\halflength
       \advance\ypos by\gap
       \put(\xpos,\ypos){\vector(0,1){\halflength}}%
   \or \advance \ypos by-\arrowlength
       \put(\xpos,\ypos){\line(0,1)\halflength}%
       \advance\ypos by\halflength
       \advance\ypos by\gap
       \put(\xpos,\ypos){\evector(0,1){\halflength}}%
   \fi
\fi
}

\def\putmorphism(#1)(#2,#3)[#4`#5`#6]#7#8#9{{%
\run #2
\rise #3
\ifnum\rise=0
  \puthmorphism(#1)[#4`#5`#6]{#7}{#8}#9%
\else\ifnum\run=0
  \putvmorphism(#1)[#4`#5`#6]{#7}{#8}#9%
\else
\setpos(#1)%
\arrowlength #7
\arrowtype #8
\ifnum\run=0
\else\ifnum\rise=0
\else
\ifnum\run>0
    \coefa=1
\else
   \coefa=-1
\fi
\ifnum\arrowtype>0
   \coefb=0
   \coefc=-1
\else
   \coefb=\coefa
   \coefc=1
   \arrowtype=-\arrowtype
\fi
\width=2
\multiply \width by\run
\divide \width by\rise
\ifnum \width<0  \width=-\width\fi
\advance\width by60
\if l#9 \width=-\width\fi
\putbox(\xpos,\ypos){#4}
{\multiply \coefa by\arrowlength
\advance\xpos by\coefa
\multiply \coefa by\rise
\divide \coefa by\run
\advance \ypos by\coefa
\putbox(\xpos,\ypos){#5} }%
{\multiply \coefa by\arrowlength
\divide \coefa by2
\advance \xpos by\coefa
\advance \xpos by\width
\multiply \coefa by\rise
\divide \coefa by\run
\advance \ypos by\coefa
\if l#9%
   \putrbox(\xpos,\ypos){#6}%
\else\if r#9%
   \putlbox(\xpos,\ypos){#6}%
\fi\fi }%
{\multiply \rise by-\coefc
\multiply \run by-\coefc
\multiply \coefb by\arrowlength
\advance \xpos by\coefb
\multiply \coefb by\rise
\divide \coefb by\run
\advance \ypos by\coefb
\multiply \coefc by70
\advance \ypos by\coefc
\multiply \coefc by\run
\divide \coefc by\rise
\advance \xpos by\coefc
\multiply \coefa by140
\multiply \coefa by\run
\divide \coefa by\rise
\advance \arrowlength by\coefa
\ifcase\arrowtype
\or \put(\xpos,\ypos){\vector(\run,\rise){\arrowlength}}%
\or \put(\xpos,\ypos){\mvector(\run,\rise){\arrowlength}}%
\or \put(\xpos,\ypos){\evector(\run,\rise){\arrowlength}}%
\fi}\fi\fi\fi\fi}}

\newcount\numbdashes \newcount\lengthdash \newcount\increment

\def\howmanydashes{
\numbdashes=\arrowlength \lengthdash=40
\divide\numbdashes by \lengthdash
\lengthdash=\arrowlength
\divide\lengthdash by \numbdashes
\increment=\lengthdash
\multiply\lengthdash by 3
\divide\lengthdash by 5
}

\def\putdashvector(#1)(#2,#3)#4#5{%
\ifnum#3=0 \putdashhvector(#1){#4}#5
\else
\ifnum#2=0
\putdashvvector(#1){#4}#5\fi\fi}

\def\putdashhvector(#1,#2)#3#4{{%
\arrowlength=#3 \howmanydashes
\multiput(#1,#2)(\increment,0){\numbdashes}%
{\vrule height .4pt width \lengthdash\unitlength}
\arrowtype=#4 \xpos=#1
\ifnum\arrowtype<0 \advance\arrowtype by 7 \fi
\ifcase\arrowtype
\or \advance\xpos by 10
    \put(\xpos,#2){\vector(-1,0){\lengthdash}}
    \advance\xpos by 40
    \put(\xpos,#2){\vector(-1,0){\lengthdash}}
\or \advance \xpos by 10
    \put(\xpos,#2){\vector(-1,0){\lengthdash}}
    \advance\xpos by  \arrowlength
    \advance\xpos by  -50
    \put(\xpos,#2){\vector(-1,0){\lengthdash}}
\or \advance\xpos by 10
    \put(\xpos,#2){\vector(-1,0){\lengthdash}}
\or \advance\xpos by \arrowlength
    \advance\xpos by -\lengthdash
    \put(\xpos,#2){\vector(1,0){\lengthdash}}
\or {\advance\xpos by 10
    \put(\xpos,#2){\vector(1,0){\lengthdash}}}
    \advance\xpos by \arrowlength
    \advance\xpos by -\lengthdash
    \put(\xpos,#2){\vector(1,0){\lengthdash}}
\or \advance\xpos by \arrowlength
    \advance\xpos by -\lengthdash
    \put(\xpos,#2){\vector(1,0){\lengthdash}}
    \advance\xpos by -40
    \put(\xpos,#2){\vector(1,0){\lengthdash}}
   \fi
}}

\def\putdashvvector(#1,#2)#3#4{{%
\arrowlength=#3 \howmanydashes
\ypos=#2 \advance\ypos by -\arrowlength
\multiput(#1,#2)(0,\increment){\numbdashes}%
    {\vrule width .4pt height \lengthdash\unitlength}
\arrowtype=#4 \ypos=#2
\ifnum\arrowtype<0 \advance\arrowtype by 7 \fi
\ifcase\arrowtype
\or \advance\ypos by \arrowlength \advance\ypos by -40
    \put(#1,\ypos){\vector(0,1){\lengthdash}}
    \advance\ypos by -40
    \put(#1,\ypos){\vector(0,1){\lengthdash}}
\or \advance\ypos by 10
    \put(#1,\ypos){\vector(0,1){\lengthdash}}
    \advance\ypos by \arrowlength \advance\ypos by -40
    \put(#1,\ypos){\vector(0,1){\lengthdash}}
\or \advance\ypos by \arrowlength \advance\ypos by -40
    \put(#1,\ypos){\vector(0,1){\lengthdash}}
\or \advance\ypos by 10
    \put(#1,\ypos){\vector(0,-1){\lengthdash}}
\or \advance\ypos by 10
    \put(#1,\ypos){\vector(0,-1){\lengthdash}}
    \advance\ypos by \arrowlength \advance\ypos by -40
    \put(#1,\ypos){\vector(0,-1){\lengthdash}}
\or \advance\ypos by 10
    \put(#1,\ypos){\vector(0,-1){\lengthdash}}
    \advance\ypos by 40
    \put(#1,\ypos){\vector(0,-1){\lengthdash}}
\fi
}}

\def\puthmorphism(#1,#2)[#3`#4`#5]#6#7#8{{%
\xpos #1
\ypos #2
\width #6
\arrowlength #6
\arrowtype=#7
\putbox(\xpos,\ypos){#3\vphantom{#4}}%
{\advance \xpos by\arrowlength
\putbox(\xpos,\ypos){\vphantom{#3}#4}}%
\horsize{\tempcounta}{#3}%
\horsize{\tempcountb}{#4}%
\divide \tempcounta by2
\divide \tempcountb by2
\advance \tempcounta by30
\advance \tempcountb by30
\advance \xpos by\tempcounta
\advance \arrowlength by-\tempcounta
\advance \arrowlength by-\tempcountb
\putvector(\xpos,\ypos)(1,0)\arrowlength\arrowtype
\divide \arrowlength by2
\advance \xpos by\arrowlength
\vertsize{\tempcounta}{#5}%
\divide\tempcounta by2
\advance \tempcounta by20
\if a#8 %
   \advance \ypos by\tempcounta
   \putbox(\xpos,\ypos){#5}%
\else
   \advance \ypos by-\tempcounta
   \putbox(\xpos,\ypos){#5}%
\fi}}

\def\putvmorphism(#1,#2)[#3`#4`#5]#6#7#8{{%
\xpos #1
\ypos #2
\arrowlength #6
\arrowtype #7
\settowidth{\xlen}{$#5$}%
\putbox(\xpos,\ypos){#3}%
{\advance \ypos by-\arrowlength
\putbox(\xpos,\ypos){#4}}%
{\advance\arrowlength by-140
\advance \ypos by-70
\ifdim\xlen>0pt
   \if m#8%
      \putsplitvector(\xpos,\ypos)\arrowlength\arrowtype
   \else
   \putvector(\xpos,\ypos)(0,-1)\arrowlength\arrowtype
   \fi
\else
   \putvector(\xpos,\ypos)(0,-1)\arrowlength\arrowtype
\fi}%
\ifdim\xlen>0pt
   \divide \arrowlength by2
   \advance\ypos by-\arrowlength
   \if l#8%
      \advance \xpos by-40
      \putrbox(\xpos,\ypos){#5}%
   \else\if r#8%
      \advance \xpos by40
      \putlbox(\xpos,\ypos){#5}%
   \else
      \putbox(\xpos,\ypos){#5}%
   \fi\fi
\fi
}}

\def\putsquarep<#1>(#2)[#3;#4`#5`#6`#7]{{%
\setsqparms[#1]%
\setpos(#2)%
\settokens[#3]%
\puthmorphism(\xpos,\ypos)[\tokenc`\tokend`{#7}]{\width}{\arrowtyped}b%
\advance\ypos by \height
\puthmorphism(\xpos,\ypos)[\tokena`\tokenb`{#4}]{\width}{\arrowtypea}a%
\putvmorphism(\xpos,\ypos)[``{#5}]{\height}{\arrowtypeb}l%
\advance\xpos by \width
\putvmorphism(\xpos,\ypos)[``{#6}]{\height}{\arrowtypec}r%
}}

\def\putsquare{\@ifnextchar <{\putsquarep}{\putsquarep%
   <\arrowtypea`\arrowtypeb`\arrowtypec`\arrowtyped;\width`\height>}}
\def\square{\@ifnextchar< {\squarep}{\squarep
   <\arrowtypea`\arrowtypeb`\arrowtypec`\arrowtyped;\width`\height>}}
\def\squarep<#1>[#2`#3`#4`#5;#6`#7`#8`#9]{{
\setsqparms[#1]
\diagram
\putsquarep<\arrowtypea`\arrowtypeb`\arrowtypec`
\arrowtyped;\width`\height>
(0,0)[#2`#3`#4`{#5};#6`#7`#8`{#9}]
\enddiagram
}}                                                 
\def\putptrianglep<#1>(#2,#3)[#4`#5`#6;#7`#8`#9]{{%
\settriparms[#1]%
\xpos=#2 \ypos=#3
\advance\ypos by \height
\puthmorphism(\xpos,\ypos)[#4`#5`{#7}]{\height}{\arrowtypea}a%
\putvmorphism(\xpos,\ypos)[`#6`{#8}]{\height}{\arrowtypeb}l%
\advance\xpos by\height
\putmorphism(\xpos,\ypos)(-1,-1)[``{#9}]{\height}{\arrowtypec}r%
}}

\def\putptriangle{\@ifnextchar <{\putptrianglep}{\putptrianglep
   <\arrowtypea`\arrowtypeb`\arrowtypec;\height>}}
\def\ptriangle{\@ifnextchar <{\ptrianglep}{\ptrianglep
   <\arrowtypea`\arrowtypeb`\arrowtypec;\height>}}
\def\ptrianglep<#1>[#2`#3`#4;#5`#6`#7]{{
\settriparms[#1]
\diagram
\putptrianglep<\arrowtypea`\arrowtypeb`
\arrowtypec;\height>
(0,0)[#2`#3`#4;#5`#6`{#7}]
\enddiagram
}}                                            

\def\putqtrianglep<#1>(#2,#3)[#4`#5`#6;#7`#8`#9]{{%
\settriparms[#1]%
\xpos=#2 \ypos=#3
\advance\ypos by\height
\puthmorphism(\xpos,\ypos)[#4`#5`{#7}]{\height}{\arrowtypea}a%
\putmorphism(\xpos,\ypos)(1,-1)[``{#8}]{\height}{\arrowtypeb}l%
\advance\xpos by\height
\putvmorphism(\xpos,\ypos)[`#6`{#9}]{\height}{\arrowtypec}r%
}}

\def\putqtriangle{\@ifnextchar <{\putqtrianglep}{\putqtrianglep
   <\arrowtypea`\arrowtypeb`\arrowtypec;\height>}}
\def\qtriangle{\@ifnextchar <{\qtrianglep}{\qtrianglep
   <\arrowtypea`\arrowtypeb`\arrowtypec;\height>}}
\def\qtrianglep<#1>[#2`#3`#4;#5`#6`#7]{{
\settriparms[#1]
\width=\height                                
\diagram
\putqtrianglep<\arrowtypea`\arrowtypeb`
\arrowtypec;\height>
(0,0)[#2`#3`#4;#5`#6`{#7}]
\enddiagram
}}

\def\putdtrianglep<#1>(#2,#3)[#4`#5`#6;#7`#8`#9]{{%
\settriparms[#1]%
\xpos=#2 \ypos=#3
\puthmorphism(\xpos,\ypos)[#5`#6`{#9}]{\height}{\arrowtypec}b%
\advance\xpos by \height \advance\ypos by\height
\putmorphism(\xpos,\ypos)(-1,-1)[``{#7}]{\height}{\arrowtypea}l%
\putvmorphism(\xpos,\ypos)[#4``{#8}]{\height}{\arrowtypeb}r%
}}

\def\putdtriangle{\@ifnextchar <{\putdtrianglep}{\putdtrianglep
   <\arrowtypea`\arrowtypeb`\arrowtypec;\height>}}
\def\dtriangle{\@ifnextchar <{\dtrianglep}{\dtrianglep
   <\arrowtypea`\arrowtypeb`\arrowtypec;\height>}}
\def\dtrianglep<#1>[#2`#3`#4;#5`#6`#7]{{
\settriparms[#1]
\width=\height                                
\diagram
\putdtrianglep<\arrowtypea`\arrowtypeb`
\arrowtypec;\height>
(0,0)[#2`#3`#4;#5`#6`{#7}]
\enddiagram
}}

\def\putbtrianglep<#1>(#2,#3)[#4`#5`#6;#7`#8`#9]{{%
\settriparms[#1]%
\xpos=#2 \ypos=#3
\puthmorphism(\xpos,\ypos)[#5`#6`{#9}]{\height}{\arrowtypec}b%
\advance\ypos by\height
\putmorphism(\xpos,\ypos)(1,-1)[``{#8}]{\height}{\arrowtypeb}r%
\putvmorphism(\xpos,\ypos)[#4``{#7}]{\height}{\arrowtypea}l%
}}

\def\putbtriangle{\@ifnextchar <{\putbtrianglep}{\putbtrianglep
   <\arrowtypea`\arrowtypeb`\arrowtypec;\height>}}
\def\btriangle{\@ifnextchar <{\btrianglep}{\btrianglep
   <\arrowtypea`\arrowtypeb`\arrowtypec;\height>}}
\def\btrianglep<#1>[#2`#3`#4;#5`#6`#7]{{
\settriparms[#1]
\width=\height                               
\diagram
\putbtrianglep<\arrowtypea`\arrowtypeb`
\arrowtypec;\height>
(0,0)[#2`#3`#4;#5`#6`{#7}]
\enddiagram
}}

\def\putAtrianglep<#1>(#2,#3)[#4`#5`#6;#7`#8`#9]{{%
\settriparms[#1]%
\xpos=#2 \ypos=#3
{\multiply \height by2
\puthmorphism(\xpos,\ypos)[#5`#6`{#9}]{\height}{\arrowtypec}b}%
\advance\xpos by\height \advance\ypos by\height
\putmorphism(\xpos,\ypos)(-1,-1)[#4``{#7}]{\height}{\arrowtypea}l%
\putmorphism(\xpos,\ypos)(1,-1)[``{#8}]{\height}{\arrowtypeb}r%
}}

\def\putAtriangle{\@ifnextchar <{\putAtrianglep}{\putAtrianglep
   <\arrowtypea`\arrowtypeb`\arrowtypec;\height>}}
\def\Atriangle{\@ifnextchar <{\Atrianglep}{\Atrianglep
   <\arrowtypea`\arrowtypeb`\arrowtypec;\height>}}
\def\Atrianglep<#1>[#2`#3`#4;#5`#6`#7]{{
\settriparms[#1]
\width=\height                                     
\diagram
\putAtrianglep<\arrowtypea`\arrowtypeb`
\arrowtypec;\height>
(0,0)[#2`#3`#4;#5`#6`{#7}]
\enddiagram
}}

\def\putAtrianglepairp<#1>(#2)[#3;#4`#5`#6`#7`#8]{{%
\settripairparms[#1]%
\setpos(#2)%
\settokens[#3]%
\puthmorphism(\xpos,\ypos)[\tokenb`\tokenc`{#7}]{\height}{\arrowtyped}b%
\advance\xpos by\height
\puthmorphism(\xpos,\ypos)[\phantom{\tokenc}`\tokend`{#8}]%
{\height}{\arrowtypee}b%
\advance\ypos by\height
\putmorphism(\xpos,\ypos)(-1,-1)[\tokena``{#4}]{\height}{\arrowtypea}l%
\putvmorphism(\xpos,\ypos)[``{#5}]{\height}{\arrowtypeb}m%
\putmorphism(\xpos,\ypos)(1,-1)[``{#6}]{\height}{\arrowtypec}r%
}}

\def\putAtrianglepair{\@ifnextchar <{\putAtrianglepairp}{\putAtrianglepairp%
   <\arrowtypea`\arrowtypeb`\arrowtypec`\arrowtyped`\arrowtypee;\height>}}
\def\Atrianglepair{\@ifnextchar <{\Atrianglepairp}{\Atrianglepairp%
   <\arrowtypea`\arrowtypeb`\arrowtypec`\arrowtyped`\arrowtypee;\height>}}

\def\Atrianglepairp<#1>[#2;#3`#4`#5`#6`#7]{{
\settripairparms[#1]
\settokens[#2]
\width=\height                                
\diagram
\putAtrianglepairp                            
<\arrowtypea`\arrowtypeb`\arrowtypec`
\arrowtyped`\arrowtypee;\height>
(0,0)[{#2};#3`#4`#5`#6`{#7}]
\enddiagram
}}

\def\putVtrianglep<#1>(#2,#3)[#4`#5`#6;#7`#8`#9]{{%
\settriparms[#1]%
\xpos=#2 \ypos=#3
\advance\ypos by\height
{\multiply\height by2
\puthmorphism(\xpos,\ypos)[#4`#5`{#7}]{\height}{\arrowtypea}a}%
\putmorphism(\xpos,\ypos)(1,-1)[`#6`{#8}]{\height}{\arrowtypeb}l%
\advance\xpos by\height
\advance\xpos by\height
\putmorphism(\xpos,\ypos)(-1,-1)[``{#9}]{\height}{\arrowtypec}r%
}}

\def\putVtriangle{\@ifnextchar <{\putVtrianglep}{\putVtrianglep
   <\arrowtypea`\arrowtypeb`\arrowtypec;\height>}}
\def\Vtriangle{\@ifnextchar <{\Vtrianglep}{\Vtrianglep
   <\arrowtypea`\arrowtypeb`\arrowtypec;\height>}}
\def\Vtrianglep<#1>[#2`#3`#4;#5`#6`#7]{{
\settriparms[#1]
\width=\height                                 
\diagram
\putVtrianglep<\arrowtypea`\arrowtypeb`
\arrowtypec;\height>
(0,0)[#2`#3`#4;#5`#6`{#7}]
\enddiagram
}}

\def\putVtrianglepairp<#1>(#2)[#3;#4`#5`#6`#7`#8]{{
\settripairparms[#1]%
\setpos(#2)%
\settokens[#3]%
\advance\ypos by\height
\putmorphism(\xpos,\ypos)(1,-1)[`\tokend`{#6}]{\height}{\arrowtypec}l%
\puthmorphism(\xpos,\ypos)[\tokena`\tokenb`{#4}]{\height}{\arrowtypea}a%
\advance\xpos by\height
\puthmorphism(\xpos,\ypos)[\phantom{\tokenb}`\tokenc`{#5}]%
{\height}{\arrowtypeb}a%
\putvmorphism(\xpos,\ypos)[``{#7}]{\height}{\arrowtyped}m%
\advance\xpos by\height
\putmorphism(\xpos,\ypos)(-1,-1)[``{#8}]{\height}{\arrowtypee}r%
}}

\def\putVtrianglepair{\@ifnextchar <{\putVtrianglepairp}{\putVtrianglepairp%
    <\arrowtypea`\arrowtypeb`\arrowtypec`\arrowtyped`\arrowtypee;\height>}}
\def\Vtrianglepair{\@ifnextchar <{\Vtrianglepairp}{\Vtrianglepairp%
    <\arrowtypea`\arrowtypeb`\arrowtypec`\arrowtyped`\arrowtypee;\height>}}
\def\Vtrianglepairp<#1>[#2;#3`#4`#5`#6`#7]{{
\settripairparms[#1]
\settokens[#2]
\diagram
\putVtrianglepairp                             
<\arrowtypea`\arrowtypeb`\arrowtypec`
\arrowtyped`\arrowtypee;\height>
(0,0)[{#2};#3`#4`#5`#6`{#7}]
\enddiagram
}}

\def\putCtrianglep<#1>(#2,#3)[#4`#5`#6;#7`#8`#9]{{%
\settriparms[#1]%
\xpos=#2 \ypos=#3
\advance\ypos by\height
\putmorphism(\xpos,\ypos)(1,-1)[``{#9}]{\height}{\arrowtypec}l%
\advance\xpos by\height
\advance\ypos by\height
\putmorphism(\xpos,\ypos)(-1,-1)[#4`#5`{#7}]{\height}{\arrowtypea}l%
{\multiply\height by 2
\putvmorphism(\xpos,\ypos)[`#6`{#8}]{\height}{\arrowtypeb}r}%
}}

\def\putCtriangle{\@ifnextchar <{\putCtrianglep}{\putCtrianglep
    <\arrowtypea`\arrowtypeb`\arrowtypec;\height>}}
\def\Ctriangle{\@ifnextchar <{\Ctrianglep}{\Ctrianglep
    <\arrowtypea`\arrowtypeb`\arrowtypec;\height>}}
\def\Ctrianglep<#1>[#2`#3`#4;#5`#6`#7]{{
\settriparms[#1]
\width=\height                               
\diagram
\putCtrianglep<\arrowtypea`\arrowtypeb`
\arrowtypec;\height>
(0,0)[#2`#3`#4;#5`#6`{#7}]
\enddiagram
}}                                           
\def\putDtrianglep<#1>(#2,#3)[#4`#5`#6;#7`#8`#9]{{%
\settriparms[#1]%
\xpos=#2 \ypos=#3
\advance\xpos by\height \advance\ypos by\height
\putmorphism(\xpos,\ypos)(-1,-1)[``{#9}]{\height}{\arrowtypec}r%
\advance\xpos by-\height \advance\ypos by\height
\putmorphism(\xpos,\ypos)(1,-1)[`#5`{#8}]{\height}{\arrowtypeb}r%
{\multiply\height by 2
\putvmorphism(\xpos,\ypos)[#4`#6`{#7}]{\height}{\arrowtypea}l}%
}}

\def\putDtriangle{\@ifnextchar <{\putDtrianglep}{\putDtrianglep
    <\arrowtypea`\arrowtypeb`\arrowtypec;\height>}}
\def\Dtriangle{\@ifnextchar <{\Dtrianglep}{\Dtrianglep
   <\arrowtypea`\arrowtypeb`\arrowtypec;\height>}}
\def\Dtrianglep<#1>[#2`#3`#4;#5`#6`#7]{{
\settriparms[#1]
\width=\height                              
\diagram
\putDtrianglep<\arrowtypea`\arrowtypeb`
\arrowtypec;\height>
(0,0)[#2`#3`#4;#5`#6`{#7}]
\enddiagram
}}                                          
\def\setrecparms[#1`#2]{\width=#1 \height=#2}%

\def\recursep<#1`#2>[#3;#4`#5`#6`#7`#8]{{%
\width=#1 \height=#2
\settokens[#3]
\settowidth{\tempdimen}{$\tokena$}
\ifdim\tempdimen=0pt
  \savebox{\tempboxa}{\hbox{$\tokenb$}}%
  \savebox{\tempboxb}{\hbox{$\tokend$}}%
  \savebox{\tempboxc}{\hbox{$#6$}}%
\else
  \savebox{\tempboxa}{\hbox{$\hbox{$\tokena$}\times\hbox{$\tokenb$}$}}%
  \savebox{\tempboxb}{\hbox{$\hbox{$\tokena$}\times\hbox{$\tokend$}$}}%
  \savebox{\tempboxc}{\hbox{$\hbox{$\tokena$}\times\hbox{$#6$}$}}%
\fi
\ypos=\height
\divide\ypos by 2
\xpos=\ypos
\advance\xpos by \width
\bfig
\putCtrianglep<-1`1`1;\ypos>(0,0)[`\tokenc`;#5`#6`{#7}]%
\puthmorphism(\ypos,0)[\tokend`\usebox{\tempboxb}`{#8}]{\width}{-1}b%
\puthmorphism(\ypos,\height)[\tokenb`\usebox{\tempboxa}`{#4}]{\width}{-1}a%
\advance\ypos by \width
\putvmorphism(\ypos,\height)[``\usebox{\tempboxc}]{\height}1r%
\efig
}}

\def\recurse{\@ifnextchar <{\recursep}{\recursep<\width`\height>}}

\def\puttwohmorphisms(#1,#2)[#3`#4;#5`#6]#7#8#9{{%
\puthmorphism(#1,#2)[#3`#4`]{#7}0a
\ypos=#2
\advance\ypos by 20
\puthmorphism(#1,\ypos)[\phantom{#3}`\phantom{#4}`#5]{#7}{#8}a
\advance\ypos by -40
\puthmorphism(#1,\ypos)[\phantom{#3}`\phantom{#4}`#6]{#7}{#9}b
}}

\def\puttwovmorphisms(#1,#2)[#3`#4;#5`#6]#7#8#9{{%
\putvmorphism(#1,#2)[#3`#4`]{#7}0a
\xpos=#1
\advance\xpos by -20
\putvmorphism(\xpos,#2)[\phantom{#3}`\phantom{#4}`#5]{#7}{#8}l
\advance\xpos by 40
\putvmorphism(\xpos,#2)[\phantom{#3}`\phantom{#4}`#6]{#7}{#9}r
}}

\def\puthcoequalizer(#1)[#2`#3`#4;#5`#6`#7]#8#9{{%
\setpos(#1)%
\puttwohmorphisms(\xpos,\ypos)[#2`#3;#5`#6]{#8}11%
\advance\xpos by #8
\puthmorphism(\xpos,\ypos)[\phantom{#3}`#4`#7]{#8}1{#9}
}}

\def\putvcoequalizer(#1)[#2`#3`#4;#5`#6`#7]#8#9{{%
\setpos(#1)%
\puttwovmorphisms(\xpos,\ypos)[#2`#3;#5`#6]{#8}11%
\advance\ypos by -#8
\putvmorphism(\xpos,\ypos)[\phantom{#3}`#4`#7]{#8}1{#9}
}}

\def\putthreehmorphisms(#1)[#2`#3;#4`#5`#6]#7(#8)#9{{%
\setpos(#1) \settypes(#8)
\if a#9 %
     \vertsize{\tempcounta}{#5}%
     \vertsize{\tempcountb}{#6}%
     \ifnum \tempcounta<\tempcountb \tempcounta=\tempcountb \fi
\else
     \vertsize{\tempcounta}{#4}%
     \vertsize{\tempcountb}{#5}%
     \ifnum \tempcounta<\tempcountb \tempcounta=\tempcountb \fi
\fi
\advance \tempcounta by 60
\puthmorphism(\xpos,\ypos)[#2`#3`#5]{#7}{\arrowtypeb}{#9}
\advance\ypos by \tempcounta
\puthmorphism(\xpos,\ypos)[\phantom{#2}`\phantom{#3}`#4]{#7}{\arrowtypea}{#9}
\advance\ypos by -\tempcounta \advance\ypos by -\tempcounta
\puthmorphism(\xpos,\ypos)[\phantom{#2}`\phantom{#3}`#6]{#7}{\arrowtypec}{#9}
}}

\def\setarrowtoks[#1`#2`#3`#4`#5`#6]{%
\def\toka{#1}
\def\tokb{#2}
\def\tokc{#3}
\def\tokd{#4}
\def\toke{#5}
\def\tokf{#6}
}
\def\hex{\@ifnextchar <{\hexp}{\hexp<1000`400>}}
\def\hexp<#1`#2>[#3`#4`#5`#6`#7`#8;#9]{%
\setarrowtoks[#9]
\yext=#2 \advance \yext by #2
\xext=#1 \advance\xext by \yext
\bfig
\putCtriangle<-1`0`1;#2>(0,0)[`#5`;\tokb``\tokd]
\xext=#1 \yext=#2 \advance \yext by #2
\putsquare<1`0`0`1;\xext`\yext>(#2,0)[#3`#4`#7`#8;\toka```\tokf]
\advance \xext by #2
\putDtriangle<0`1`-1;#2>(\xext,0)[`#6`;`\tokc`\toke]
\efig
}

\newtheorem{theorem}{Theorem}
\newtheorem{proposition}{Proposition}
\newtheorem{lemma}{Lemma}
\newtheorem{example}{Example}
\newtheorem{remark}{Remark}

\def\mr#1{\smash{
        \mathop{\longrightarrow}\limits^{#1}}}
\def\ml#1{\smash{
        \mathop{\longleftarrow}\limits^{#1}}}
\def\md#1{\Big\downarrow
      \rlap{$\vcenter{\hbox{$\scriptstyle#1$}}$}}
\def\mdi#1{\bigcup
      \rlap{$\vcenter{\hbox{$\scriptstyle#1$}}$}}
\def\mri#1{\smash{
        \mathop{\supset}\limits^{#1}}}
\def\mli#1{\smash{
        \mathop{\subset}\limits^{#1}}}

\QQQ{Language}{
American English
}

\begin{document}

\author{S.A. Lysenko}
\title{On the functional equation $f(p(z))=g(q(z))$, where $p,q$ are
``generalized'' polynomials and $f,g$ are meromorphic functions}
\maketitle

\begin{abstract}
We find all the solutions to the equation $f(p(z))=g(q(z))$, where $p,q$ are
polynomials and $f,g$ are transcendental meromorphic functions in ${\bf C}$.
In fact, a more general algebraic problem is solved.
\end{abstract}

\tableofcontents

\sloppy

\section*{Introduction}

\addcontentsline{toc}{section}{Introduction}%
\markboth{On the functional equation $f(p(z))=g(q(z))$}{On the functional equation $f(p(z))=g(q(z))$}

\subsection{Motivation}

This paper and the previous one \cite{c3} were partially motivated by the
following result by L. Flatto \cite{c1} :

Let $p,q\in {\bf C[}z]$ be polynomials of equal degree. Let 
\begin{equation}
\label{e0.1}f\circ p=g\circ q 
\end{equation}
where $f$ and $g$ are nonconstant entire functions on ${\bf C}$. Then one of
the following is true :

\begin{description}
\item[(i)]  $p\left( z\right) =\lambda q\left( z\right) +a$, with $\lambda
,a\in {\bf C}$;

\item[(ii)]  $p\left( z\right) =r\left( z\right) ^2+a$,\ $q\left( z\right)
=br\left( z\right) ^2+cr\left( z\right) +d$, where $r$ is a polynomial in $z$
and $a,b,c,d,\in {\bf C}$,\ $b\neq 0$.
\end{description}

L. Flatto asked \cite[question 5]{c4} whether there is an analog of his
theorem if $\deg p\neq \deg q.$ One can also ask what happens if $f$ and $g$
are not entire but meromorphic functions (on the whole ${\bf C}$ or only in
neighborhood of infinity). Partial results related to Flatto's question were
obtained in \cite{c11}, \cite{c2}, \cite{c12}, \cite{c3} (\cite{c2} contains
a survey of most of these results). The goal of this paper is to describe
all pairs $p,q$ for which there exist nonconstant meromorphic $f$ and $g$
satisfying (\ref{e0.1}) and there exist no rational $f$ and $g$ with this
property (actually we consider a more general problem; see \ref{ch0.2} and
Section \ref{Sec1}).

Our interest to equation (\ref{e0.1}) is also motivated by its relation to
the following problem which seems interesting: describe equivalence
relations ${\bf R}$ on ${\bf C}$ such that 1) ${\bf R}$ considered as a
subset of ${\bf C}^2$ is a union of a sequence of algebraic curves, 2) there
exists a nonconstant meromorphic function on ${\bf C}$ whose restriction to
each equivalence class of ${\bf R}$ is constant. Such equivalence relations
can be considered as generalizations of discrete subgroups of biholomorphic
automorphisms of ${\bf C}$ (a discrete subgroup $\Gamma $ defines the
following equivalence relation : $z\sim w$ iff $z=\gamma (w)$ for some $%
\gamma \in \Gamma $; clearly, this equivalence relation satisfies conditions
1) and 2)). Notice that a solution to (\ref{e0.1}) with polynomial $p,q$ and
meromorphic $f,g$ gives rise to an equivalence relation $R_{p,q}$ satisfying
1) and 2) ($R_{p,q}$ is the equivalence relation generated by $R_p=\{\left(
z,u\right) \in {\bf C}^2\mid p\left( z\right) =p\left( u\right) \}$ and $%
R_q=\{(z,u)\in {\bf C}^2\mid q(z)=q(u)\}$). In some sense $R_{p,q}$ is the
equivalence relation associated to the group $\Gamma $ generated by $%
z\mapsto h_p(z):=p^{-1}(p(z))$ and $z\mapsto h_q(z):=q^{-1}(q(z))$. Notice
however that $h_p$ and $h_q$ are holomorphic only in a neighborhood of $%
\infty $, so $\Gamma $ is a group of ${\em germs}$\/\ of conformal mappings
: $({\bf \bar C},\infty )\rightarrow ({\bf \bar C},\infty )$.

\subsection{Formulation of the problem\label{ch0.2}}

Let $\left( X,\infty _X\right) $, $\left( Y,\infty _Y\right) $ be compact
Riemann surfaces with marked points $\infty _X\in X,\;\infty _Y\in Y$. By
abuse of notation we write $\infty $ instead of $\infty _X$ and $\infty _Y$.
Let $p:\left( X,\infty \right) \rightarrow \left( Y,\infty \right) $ be a
holomorphic map. We say that $p$ is a {\em generalized polynomial}${\em \/}$
if $p^{-1}\left( \left\{ \infty \right\} \right) =\left\{ \infty \right\} $.

We will study equation (\ref{e0.1}), where $p:\left( X,\infty _X\right)
\rightarrow \left( Y,\infty _Y\right) $ and $q:\left( X,\infty _X\right)
\rightarrow \left( Z,\infty _Z\right) $ are generalized polynomials and $f,g$
are meromorphic functions in punctured neighborhoods of $\infty _Y$ and $%
\infty _Z$ respectively. By rational function on a compact Riemann surface
we shall mean a meromorphic function on it (this agrees with the usage of
the term ``rational function'' in algebraic geometry). It is required to
find all pairs $p,q$ such that there exists nonconstant $f,g$ satisfying (%
\ref{e0.1}) and there exist no rational $f,g$ with this property. In fact,
we solve a more general algebraic problem, which is explained in Section \ref
{Sec1}.

\subsection{Main result}

There are several standard solutions to (\ref{e0.1}).

\begin{example}
\label{x0.1}Let $p\left( z\right) =z^n$,\ $q\left( z\right) =\left(
z+1\right) ^m$ with $n,m,\limfunc{lcm}\left( n,m\right) \in \left\{
2,3,4,6\right\} $. Then there exist nonconstant functions $f,g$ meromorphic
in ${\bf C}$ and satisfying (\ref{e0.1}). There exist no rational $f,g$ with
this property.
\end{example}

\begin{remark}
\label{r0.1}Suppose we are given a diagram of generalized polynomials%
$$
\Vtrianglep<0`1`1;400>[(X,\infty )`(Z,\infty )`(Y,\infty );`p`q] 
$$
such that $\limfunc{gcd}(\deg p,\deg q)=1$. Then there exists a diagram of
generalized polynomials%
$$
\TeXButton{Diagram p.2 n.2}
{
    \matrix{{(W,\infty)}&\mr{p_1}&{(Z,\infty)}\cr
      \md{q_1}&&\md{q}\cr
      {(X,\infty)}&\mr{p}&{(Y,\infty)}\cr}
} 
$$
such that $\deg p_1=\deg p$, $\deg q_1=\deg q$. It is unique up to canonical
isomorphism.In fact, $W=W_0$, where $W_0$ is the normalization
(=desingularization) of the analytic curve $X\times _YZ=\{\left( x,z\right)
\in X\times Z\mid p\left( x\right) =q\left( z\right) \}$ (Let us explain
that if $\deg p$ and $\deg q$ are coprime then $W_0$ has only one point over 
$\infty _Y$, which implies that $W_0$ is connected; so the maps $%
W_0\rightarrow X$ and $W_0\rightarrow Y$ are generalized polynomials).
Notice that if $X,Y,Z$ are of genus $0$ then, as a rule, $W$ is of genus
greater than $0$.
\end{remark}

\begin{example}
\label{x0.2}Let $\tilde p,\tilde q$ be the pair of polynomials from Example 
\ref{x0.1}, $\deg \tilde p=n$, $\deg \tilde q=m$. Let $h:\left( Y,\infty
\right) \rightarrow ({\bf CP}^1,\infty )$, $r:\left( Z,\infty \right)
\rightarrow ({\bf CP}^1,\infty )$ be generalized polynomials, $\deg h=\alpha 
$, $\deg r=\beta $. Suppose that $\limfunc{gcd}\left( \alpha ,n\right) =%
\limfunc{gcd}(\beta ,m)=\limfunc{gcd}(\alpha ,\beta )=1$. Using Remark \ref
{r0.1}, we get the following commutative diagram of generalized polynomials%
$$
\TeXButton{Diagram p.3 n.1}
{
    \matrix{W&\mr{h_2}&Z_1&\mr{q_1}&Z\cr
               \md{r_2}&&\md{r_1}&&\md{r}\cr
               Y_1&\mr{h_1}&\bf{CP}^1&\mr{\tilde q}&\bf{CP}^1\cr
               \md{p_1}&&\md{\tilde p}&&&\cr
               Y&\mr{h}&\bf{CP}^1&&&\cr}
} 
$$
with $\deg r_1=\deg r_2=\beta $, $\deg h_1=\deg h_2=\alpha $, $\deg p_1=n$, $%
\deg q_1=m$. Let $s:\left( X,\infty \right) \rightarrow \left( W,\infty
\right) $ be a generalized polynomial. Put $p=p_1\circ r_2\circ s$, $%
q=q_1\circ h_2\circ s$. The pair $p,q$ is a solution of the problem under
consideration.
\end{example}

Our main result is that {\em Example \ref{x0.2} provides all the solutions
of our problem} (in the special case $\deg p=\deg q$ it means that $%
h,h_1,h_2 $ and $r,r_1,r_2$ are isomorphisms, so actually $p=\tilde p\circ s$%
, $q=\tilde q\circ s$).

\subsection{Organization}

In Section \ref{Sec1} we replace our problem by a more general algebraic
one, which is actually treated, and formulate the corresponding results. We
introduce the concept of irreducible pair of generalized polynomials and
separate the results in two parts. First, we reduce our problem to that with
irreducible pairs of generalized polynomials (Section \ref{Sec3}). Secondly,
we study the irreducible pairs (Section \ref{Sec4}). In Section \ref{Sec4}
we formulate the Main group-theoretic lemma, which plays a central role in
the proof of our main result. This lemma is proved in Section \ref{Sec5}.

\subsection{Conventions}

All Riemann surfaces are supposed to be connected. Recall that the following
three concepts are equivalent : a compact Riemann surface, a nonsingular
connected projective algebraic curve over ${\bf C}$, a finitely generated
field over ${\bf C}$ of transcendence degree $1$. We shall identify a point
of a curve and the corresponding place of the field of rational functions on
this curve.

\subsection{Acknowledgments}

I am greatly indebted to V.G. Drinfeld for his valuable advises and constant
attention to this work.

The author was partially supported by INTAS grant No. 94--4720.

\section{Formulation of results\label{Sec1}}

Denote by $J$ the group of (all) germs of conformal mappings : $({\bf CP}%
^1,\infty )\rightarrow ({\bf CP}^1,\infty )$.

\begin{description}
\item[Definition.]  Suppose $\Gamma $ is a subgroup of $J$. We say that $%
\Gamma $ is {\em discrete} if there exists a nonconstant function $F$
meromorphic in a punctured neighborhood of infinity in ${\bf C}$ such that $%
F(g(z))=F(z)$ for all $g\in \Gamma $.
\end{description}

In the paper \cite{c3} a necessary condition for a group $\Gamma $ to be
discrete was obtained using the results from analytic local dynamics \cite
{c13}. This condition will serve as the main tool in the proofs of our
theorems. Let us formulate it here. For any $g\in {\bf C((}\frac 1z))$ we
write $\limfunc{ord}_\infty g=n$ if $g=\sum\limits_{k=n}^\infty a_kz^{-k}$, $%
a_n\neq 0$. Put $J_k=\left\{ g\in J\mid \limfunc{ord}_\infty (g(z)-z)\geq
1-k\right\} $ for $k\leq 1$. We have $J\supset J_1\supset J_0\supset
J_{-1}\supset \ldots $ Here $J_k$ is a normal subgroup of $J$. If $\Gamma
\subset J$ is a subgroup, then we write $\Gamma _k=\Gamma \cap J_k$, $k\leq
1 $.

\begin{theorem}[\cite{c3}]
\label{th1.1}Suppose $\Gamma \subset J$ is a discrete subgroup; then

\begin{enumerate}
\item  \label{it1} at most one of the quotient groups $\Gamma _k/\Gamma
_{k-1}$ ($k\leq 1$) is nontrivial,

\item  \label{it2} for all $k\leq 1$ the subgroup $\Gamma _k/\Gamma
_{k-1}\subset J_k/J_{k-1}\simeq ({\bf C},+)$ is discrete.
\end{enumerate}
\end{theorem}

\begin{description}
\item[Definition.]  A subgroup $\Gamma \subset J$ is {\em formally discrete}
if it satisfies the conditions \ref{it1}), \ref{it2}) from Theorem \ref
{th1.1}.

\item[Remarks.] 
\begin{enumerate}
\item  The quotient group $\Gamma /\Gamma _1$ is ignored here.

\item  Theorem \ref{th1.1} can be partially proved using the result of
Scherbakov \cite{c14}.
\end{enumerate}
\end{description}

Suppose $X$ is a Riemann surface, $\infty $ is a point of $X$. Denote by $%
J(X,\infty )$ the group of germs of conformal mappings : $(X,\infty
)\rightarrow (X,\infty )$. One can identify $J(X,\infty )$ with $J$ by
choosing a local parameter at $\infty $ on $X$. Let $Y$ be another Riemann
surface and $f$ a holomorphic map from a punctured neighborhood of infinity
in $X$ to $Y$. Then we define a group $T_f$ by the formula $T_f=\left\{ g\in
J(X,\infty )\mid f\circ g=f\right\} $.

Suppose we are given compact Riemann surfaces $X,Y$ and two points $\infty
\in X$, $\infty \in Y$. Recall that a {\em generalized polynomial }\ is a
holomorphic map $p:(X,\infty )\rightarrow (Y,\infty )$ such that $%
p^{-1}(\{\infty \})=\{\infty \}$. It is easy to see that if $p$ is a
generalized polynomial then $T_p$ is a cyclic group of order $\deg p$.

Consider the following diagram of generalized polynomials : 
\begin{equation}
\label{e1.1}\TeXButton{Diagram 1.1 p.2 n.1}
{\matrix{\Atrianglep<1`1`0;400>[(X,\infty )`(Y,\infty )`(Z,\infty );p`q`]}} 
\end{equation}

Let $f$ and $g$ be meromorphic functions in punctured neighborhoods of
infinity in $Y$ and $Z$ respectively. Suppose we have 
\begin{equation}
\label{e1.2}f\circ p=g\circ q 
\end{equation}

Assume that $f$ and $g$ are nonconstant. Obviously, then $T_p$ and $T_q$
generate a discrete subgroup in $J(X,\infty )$. Conversely, if $T_p$ and $%
T_q $ generate a discrete group, then there exist nonconstant functions $f$
and $g$ as above such that (\ref{e1.2}) holds.

In fact, in some sense, we find all pairs of generalized polynomials $p$ and 
$q$ such that $T_p$ and $T_q$ generate a {\em formally} discrete group.

\begin{remark}
\label{r1.2}Denote by $\bar J$ the group of all formal diffeomorphisms : $(%
{\bf CP}^1,\infty )\rightarrow ({\bf CP}^1,\infty )$, i.e., $\bar J=\left\{
z\mapsto a_1z+a_0+a_{-1}z^{-1}+\ldots \mid a_i\in {\bf C},a_1\neq 0\right\} $
with respect to superposition. The subgroups $\bar J_k\subset \bar J$ are
defined in the same way as $J_k\subset J$. We have $J\subset \bar J$. If $%
\Gamma \subset \bar J$ is a subgroup such that $\Gamma \not \subset J$, then
the discreteness property does not make sense for $\Gamma $, whereas the
formal discreteness does. In this sense the formal discreteness is an
algebraic property.
\end{remark}

In what follows we denote by ${\cal M}(X)$ the field of meromorphic
functions on a Riemann surface $X$.

Let $p:(X,\infty )\rightarrow (Y,\infty )$ be a generalized polynomial; then
we have ${\cal M}(Y)\subset {\cal M}(X)$. Let $F$ be a field such that $%
{\cal M}(Y)\subset F\subset {\cal M}(X)$, $W$ its model, i.e., the compact
Riemann surface such that $F$ is isomorphic to ${\cal M}(W)$ over ${\bf C}$.
We get a commutative diagram :%
$$
\Vtrianglep<1`1`-1;400>[X`Y`W;p`p_1`p_2] 
$$

Put $\infty _W=p_1(\infty _X)$. Then $p_1:(X,\infty )\rightarrow (W,\infty )$
and $p_2:(W,\infty )\rightarrow (Y,\infty )$ are generalized polynomials.

The following theorem may be considered as a description of rational
solutions to the functional equation (\ref{e1.2}).

\begin{theorem}
\label{th1.2}Suppose we are given a diagram (\ref{e1.1}). Then there is an
alternative~:

\begin{enumerate}
\item  There exists a commutative diagram of generalized polynomials 
\begin{equation}
\label{e1.3}\TeXButton{Diagam 1.3 p.3 n.2}
{\matrix{\xext=1200 \yext=1250
\adjust[(X)`;`;`;(W)`]
\begin{picture}(\xext,\yext)(\xoff,\yoff)
\putAtrianglepair<1`0`1`0`0;600>(0,600)[(X,\infty)`(Y,\infty)`(V,\infty)`(Z,\infty);p``q``]
\putVtrianglepair<-1`1`1`0`1;600>(0,0)[\phantom{(Y,\infty)}`\phantom{(V,\infty)}`\phantom{(Z,\infty)}`(W,\infty);p_1`q_1```]
\putmorphism(600,1200)(0,-1)[``f]{600}1l
\putmorphism(600,600)(0,-1)[``g]{600}1l
\end{picture}
}}
\end{equation}
such that $\deg f=\limfunc{gcd}(\deg p,\deg q)$, $\deg g=(\deg p_1)\cdot
(\deg q_1)$. Such a diagram is unique up to isomorphism. The groups $T_p$
and $T_q$ generate $T_{g\circ f}$, ${\cal M}(Y)\cap {\cal M}(Z)={\cal M}(W)$%
, ${\cal M}(V)$ is the composite of the fields ${\cal M}(Y)$ and ${\cal M}(Z)
$.

\item  $T_p$ and $T_q$ generate an infinite nonabelian subgroup of $%
J(X,\infty )$. In this case ${\cal M}(Y)\cap {\cal M}(Z)={\bf C}$.
\end{enumerate}
\end{theorem}

The following theorems describe the pairs of generalized polynomials $p,q$
for which $T_p$ and $T_q$ generate an infinite formally discrete subgroup of 
$J(X,\infty )$.

\begin{proposition}
\label{pp1.1}Suppose we are given a diagram (\ref{e1.1}). Then there exists
a unique field $F_{q,p}$ with the following property. First, ${\cal M}%
(Z)\subset F_{q,p}\subset {\cal M}(X)$, $F_{q,p}\cap {\cal M}(Y)\neq {\bf C}$%
. Secondly, given a field $F$ such that ${\cal M}(Z)\subset F\subset {\cal M}%
(X)$ and $F\cap {\cal M}(Y)\neq {\bf C}$, we have $F_{q,p}\subset F$.
\end{proposition}

\begin{remark}
\label{r1.3}From the geometrical point of view it means that there exists a
commutative diagram of generalized polynomials :%
$$
\TeXButton{Diagram p.4 n.1}
{
    \matrix{X&\mr{q_1}&X_{q,p}&\mr{q_2}&Z\cr
               \md{}&&\md{}&&&\cr
               Y&\mr{}&Y_{q,p}&&&\cr}
} 
$$
such that $q_2\circ q_1=q$ and the following universal property holds. Given
a commutative diagram of generalized polynomials :%
$$
\TeXButton{Diagram p.4 n.2}
{
    \matrix{X&\mr{g}&X^{\prime}&\mr{h}&Z\cr
               \md{p}&&\md{}&&&\cr
               Y&\mr{}&Y^{\prime}&&&\cr}
} 
$$
such that $h\circ g=q$, there exists a unique holomorphic $f:X^{\prime
}\rightarrow X_{q,p}$ such that $f\circ g=q_1$.
\end{remark}

\begin{description}
\item[Definition.]  We say that the pair $p,q$ in diagram (\ref{e1.1}) is 
{\em irreducible }\ if $F_{q,p}=F_{p,q}={\cal M}(X)$.
\end{description}

\begin{remark}
\label{r1.4}Suppose the pair $p,q$ in (\ref{e1.1}) is irreducible and $\deg
p,\deg q>1$; then ${\cal M}(Y)\cap {\cal M}(Z)={\bf C}$.
\end{remark}

\begin{description}
\item[Example.]  Put $X={\bf CP}^1$, $p(z)=z^n$, $q(z)=(z+1)^m$, where $n,m$
are positive integers. Then the pair of polynomials $p,q$ is irreducible.
\end{description}

To each pair of generalized polynomials (\ref{e1.1}) we assign an
irreducible pair as follows. Put $F=F_{p,q}\cap F_{q,p}$, $F_1={\cal M}%
(Y)\cap F$, $F_2={\cal M}(Z)\cap F$. Let $K$ be the composite of $F_{p,q}$
and $F_{q,p}$. The diagram of fields commutes :%
$$
\TeXButton{Diagram p.4 n.3}
{
    \matrix{{\cal M}(X)&\mri{}&K&\mri{}&F_{q,p}&\mri{}&{\cal M}(Z)\cr
               &&\mdi{}&&\mdi{}&&\mdi{}\cr
               &&F_{p,q}&\mri{}&F&\mri{}&F_2\cr
               &&\mdi{}&&\mdi{}&&&\cr
               &&{\cal M}(Y)&\mri{}&F_1&&&\cr}
} 
$$

To this diagram there corresponds the following diagram of generalized
polynomials : 
\begin{equation}
\label{e1.4}\TeXButton{Diagram 1.4 p.4 n.4}
{
    \matrix{X&\mr{}&V&\mr{h_2}&X_{q,p}&\mr{q_1}&Z\cr
               &&\md{r_2}&&\md{r_1}&&\md{r}\cr
               &&X_{p,q}&\mr{h_1}&W&\mr{\tilde q}&W_2\cr
               &&\md{p_1}&&\md{\tilde p}&&&\cr
               &&Y&\mr{h}&W_1&&&\cr}
} 
\end{equation}

Diagram (\ref{e1.4}) will be referred to as {\em the canonical diagram}.

From the definition of $F_{p,q}$ (see Proposition \ref{pp1.1}) it follows
that $F_{p,q}$ is the composite of $F$ and ${\cal M}(Y)$. Similarly, $%
F_{q,p} $ is the composite of $F$ and ${\cal M}(Z).$ From Theorem \ref{th1.2}
it follows that $\deg h=\deg h_1=\deg h_2$, $\deg r=\deg r_1=\deg r_2$, $%
\deg p_1=\deg \tilde p$, $\deg q_1=\deg \tilde q$, $\limfunc{gcd}(\deg
h,\deg \tilde p)=\limfunc{gcd}(\deg r,\deg \tilde q)=\limfunc{gcd}(\deg
h,\deg r)=1$.

\begin{proposition}
\label{pp1.2}The pair $\tilde p,\tilde q$ is irreducible.
\end{proposition}

\begin{proposition}
\label{pp1.3}Suppose we are given a diagram (\ref{e1.1}) such that ${\cal M}%
(Y)\cap {\cal M}(Z)={\bf C}$. Let $\tilde p,\tilde q$ be the corresponding
irreducible pair of generalized polynomials. Then the following conditions
are equivalent :

\begin{itemize}
\item  $T_p$ and $T_q$ generate a discrete subgroup of $J(X,\infty )$

\item  $T_{\tilde p}$ and $T_{\tilde q}$ generate a discrete subgroup of $%
J(W,\infty )$
\end{itemize}
\end{proposition}

The assertion remains valid if we replace discreteness by formal
discreteness.

\begin{theorem}
\label{th1.3}Suppose we are given a diagram (\ref{e1.1}) such that the pair $%
p,q$ is irreducible and $\deg p,\deg q>1$. Suppose that $T_p$ and $T_q$
generate a formally discrete group. Then there exist a commutative diagram :%
$$
\TeXButton{Diagram p.5 n.1}
{
\matrix{(Y,\infty)&\ml{p}&(X,\infty)&\mr{q}&(Z,\infty)\cr
           \md{}&&\md{}&&\md{}\cr
           (\bf{CP}^1,\infty)&\ml{p_1}&(\bf{CP}^1,\infty)&\mr{q_1}&(\bf{CP}^1,\infty)\cr}}
$$
where the vertical arrows are isomorphisms and $p_1,q_1$ is the following
standard pair of polynomials : $p_1(z)=z^n$, $q_1(z)=(z+1)^m$ with $n,m,%
\limfunc{lcm}(n,m)\in \{2,3,4,6\}$. Conversely, the pair $p_1,q_1$ is
irreducible, $T_{p_1}$ and $T_{q_1}$ generate a discrete subgroup of $J$, $%
{\bf C}(p_1)\cap {\bf C}(q_1)={\bf C}$.
\end{theorem}

Our main result, which was formulated in the {\it Introduction}, follows
from Theorems \ref{th1.1}, \ref{th1.2}, \ref{th1.3} and Propositions \ref
{pp1.1}, \ref{pp1.2}, \ref{pp1.3}.

\section{An algebraic set-up\label{Sec2}}

Suppose $X$ is a compact Riemann surface, $\infty $ is a point of $X$. The
place of ${\cal M}(X)$ corresponding to $\infty $ will be denoted by the
same symbol $\infty $. We denote by ${\cal M}(X)_\infty $ the completion of $%
{\cal M}(X)$ at $\infty $. Let $p:(X,\infty )\rightarrow (Y,\infty )$ be a
generalized polynomial. Denote the restriction of $\infty $ to ${\cal M}(Y)$
by the same letter. It is known that ${\cal M}(X)_\infty $ is a cyclic
Galois extension of ${\cal M}(Y)_\infty $ of order $\deg p$. To each $g\in
J(X,\infty )$ assign the automorphism of ${\cal M}(X)_\infty $ given by $%
(gf)(x)=f(g^{-1}x)$, $f\in {\cal M}(X)_\infty $. We get an embedding of $%
J(X,\infty )$ into the group of automorphisms of the topological field $%
{\cal M}(X)_\infty $ over ${\bf C}$. In what follows $J(X,\infty )$ will be
considered as a subgroup of the latter group. This embedding induces also a
canonical isomorphism between $T_p$ and $Gal({\cal M}(X)_\infty /{\cal M}%
(Y)_\infty )$. These two groups will be identified as well.

\begin{lemma}
\label{le2.1}Let $X$ and $Y$ be compact Riemann surfaces, $f:X\rightarrow Y$
a holomorphic $n$-sheeted covering. Let $g:W\rightarrow Y$ be the least
Galois covering that can be factorized as follows :%
$$
\qtriangle<1`1`1;400>[W`X`Y;`g`f] 
$$
Let $y_0\in Y$, $f^{-1}(y_0)=\{x_1,\ldots ,x_k\}\subset X$. Suppose the
multiplicity of $f$ at $x_i$ is equal to $l_i$, $w_0\in W$, $g(w_0)=y_0$;
then the multiplicity of $g$ at $w_0$ is equal to $\limfunc{lcm}(l_1,\ldots
,l_k)$.
\end{lemma}

The following explicit construction of $W$ is useful to prove this lemma.
Let $A\subset Y$ be the set of critical values of $f$. Put $Y^{\prime
}=Y\setminus A$, $X^{\prime }=X\setminus f^{-1}(A)$. Denote by $Z^{\prime }$
the set of pairs $(y,\varphi )$, where $y\in Y^{\prime }$ and $\varphi $ is
a bijection : $f^{-1}(y)\rightarrow \{1,\ldots ,n\}$. Let $g$ be the map
from $Z^{\prime }$ to $Y^{\prime }$ such that $g(y,\varphi )=y$. A Riemann
surface structure on $Z^{\prime }$ is defined in the natural way. The group $%
S_n$ acts on $Z^{\prime }$ by biholomorphic transformations and this action
is transitive on the fibres of $g$. Let $Z$ be the smooth compactification
of $Z^{\prime }$ and $W$ a connected component of $Z$. Then $g:W\rightarrow
Y $ is the desired Galois covering.

The rest of the proof is omitted.

\begin{description}
\item[Corollary.]  Let $p:(X,\infty )\rightarrow (Y,\infty )$ be a
generalized polynomial, $K$ a finite Galois extension field of ${\cal M}(X)$
such that $K$ is not ramified over $\infty \in X$. Let $L$ be the least
Galois extension of ${\cal M}(Y)$ such that $K\subset L$. Then $L$ is not
ramified over $\infty \in X$.
\end{description}

Now we consider a diagram (\ref{e1.1}) of generalized polynomials. Fix an
algebraic closure $\overline{{\cal M}(X)}$. We construct a tower of fields $%
k_p^m,k_q^m\subset \overline{{\cal M}(X)}$, $m\geq 0$ as follows. Put $%
k_p^0=k_q^0={\cal M}(X)$. Let $k_p^m$ be the least Galois extension of $%
{\cal M}(Y)$ containing $k_q^{m-1}$. Let $k_q^m$ be the least Galois
extension of ${\cal M}(Z)$ containing $k_p^{m-1}$. One proves by induction
that $k_p^m\supset k_p^{m-1}$, $k_q^m\supset k_q^{m-1}$. By definition, $%
k_p^m\supset k_q^{m-1}$, $k_q^m\supset k_p^{m-1}$. Put $E=\bigcup%
\limits_mk_p^m=\bigcup\limits_mk_q^m$. $E$ is a field containing ${\cal M}%
(X) $ and normal over both ${\cal M}(Y)$ and ${\cal M}(Z)$. Actually $E$ is
the smallest subfield of $\overline{{\cal M}(X)}$ with this property (if $%
{\cal M}(X)\subset E^{\prime }\subset \overline{{\cal M}(X)}$ and $E^{\prime
}$ is normal over both ${\cal M}(Y)$ and ${\cal M}(Z)$ then one shows by
induction that $E^{\prime }\supset k_p^m$ and $E^{\prime }\supset k_q^m$ for
all $m$).

\begin{lemma}
\label{le2.2}For all $m\geq 0$ the fields $k_p^m$ and $k_q^m$ are not
ramified over $\infty \in X$. Therefore $E$ is not ramified over $\infty \in
X$.
\end{lemma}

\TeXButton{Proof}{\proof}This follows immediately from the previous
corollary.\TeXButton{End Proof}{\endproof}

Fix a place $\infty ^{\prime }$ of $E$ over $\infty \in X$. The choice of $%
\infty ^{\prime }$ provides an embedding $E\hookrightarrow {\cal M}%
(X)_\infty =E_{\infty ^{\prime }}$ over ${\cal M}(X)$. Since $E$ is normal
over ${\cal M}(X)$, the image of $E$ in ${\cal M}(X)_\infty $ does not
depend on the choice of $\infty ^{\prime }$.

Put $G_p=Gal(E/{\cal M}(Y))$, $G_q=Gal(E/{\cal M}(Z))$, $U=Gal(E/{\cal M}%
(X)) $. Let $G$ be the subgroup of $\limfunc{Aut}E$ generated by $G_p$ and $%
G_q$. It is well known that for any place $\omega $ of ${\cal M}(X)$
(trivial on ${\bf C}$) the action of $U$ on the set of places of $E$ over $%
\omega $ is transitive.

Denote the set of places of $E$ over $\infty $ by $S$. Clearly, $S$ is
invariant with respect to $G_p$ and $G_q$. So $G$ acts on $S$. The action of 
$U$ on $S$ is free because $E$ is not ramified over $\infty $. As explained
above this action is transitive.

It is well known that assotiating to $\sigma \in Gal({\cal M}(X)_\infty /%
{\cal M}(Y)_\infty )$ its restriction to $E\subset E_{\infty ^{\prime }}=%
{\cal M}(X)_\infty $ one obtains an isomorphism :%
$$
T_p=Gal({\cal M}(X)_\infty /{\cal M}(Y)_\infty )\tilde \rightarrow \{g\in
G_p\mid g\infty ^{\prime }=\infty ^{\prime }\} 
$$

Similarly, we have an isomorphism :%
$$
T_q=Gal({\cal M}(X)_\infty /{\cal M}(Z)_\infty )\tilde \rightarrow \{g\in
G_q\mid g\infty ^{\prime }=\infty ^{\prime }\} 
$$

Let $\Gamma $ be the subgroup of $\limfunc{Aut}{\cal M}(X)_\infty $
generated by $T_p$ and $T_q$.

\begin{lemma}
\label{le2.3}

\begin{enumerate}
\item  \label{le23.1}For any $g\in G_p$ (resp. $g\in G_q$) there exist
unique $h\in T_p$ (resp. $h\in T_q$), $\sigma \in U$ such that $g=h\sigma $.

\item  \label{le23.2}The restriction to $E$ induces an isomorphism :%
$$
\Gamma \tilde \rightarrow \{g\in G\mid g\infty ^{\prime }=\infty ^{\prime
}\} 
$$

\item  \label{le23.3}For any $g\in G$ there exist unique $h\in \Gamma $, $%
\sigma \in U$ such that $g=h\sigma $.
\end{enumerate}
\end{lemma}

\TeXButton{Proof}{\proof}

\ref{le23.1}). It follows from the fact that $S$ is identified with $G_p/T_p$
or $G_q/T_q$ and the action of $U$ on $S$ is free and transitive.

\ref{le23.2}) and \ref{le23.3}). Clearly, we have a homomorphism $f:\Gamma
\rightarrow \{g\in G\mid g\infty ^{\prime }=\infty ^{\prime }\}$. Since $E$
is dense in ${\cal M}(X)_\infty $, it follows that $f$ is injective.

From \ref{le23.1}) it follows that for every $g\in G$ there exist $h\in
\Gamma $, $\sigma \in U$ such that $g=f(h)\sigma $. These $h,\sigma $ are
unique because the action of $U$ on $S$ is free. If $g\infty ^{\prime
}=\infty ^{\prime }$ then $\sigma =1$, so $f$ is surjective.%
\TeXButton{End Proof}{\endproof}

\begin{description}
\item[Remarks.] 
\begin{enumerate}
\item  We shall consider the groups $T_p,T_q$ and $\Gamma $ as subgroups of $%
G$. At the same time $\Gamma $ can be considered as the subgroup of $%
J(X,\infty )$ generated by $T_p$ and $T_q$. By Lemma \ref{le2.3} we get a
bijection $G/U\leftrightarrow \Gamma $. The group $U$ acts on $G/U$ by left
translations, so $U$ acts on the set $\Gamma $ without preserving the group
structure of $\Gamma $. This action plays a critical role in this paper. It
has the following analytical meaning. The elements of $\Gamma $ can be
regarded as germs of algebraic functions at $\infty \in X$. The analytic
continuation around closed paths provides the monodromy action of $H$ on $%
\Gamma $, where $H$ is inverse limit of $\pi _1(X\backslash S,\infty )$, $%
S\subset X\backslash \{\infty \}$, $\#S<\infty $. There is a canonical
homomorphism $f:H\rightarrow U$ with dense image, and the monodromy action
of $H$ on $\Gamma $ comes from the action of $U$ on $\Gamma $. So the action
of $U$ on $\Gamma $ is an algebraic version of the monodromy action used in
\S 4 from \cite{c3}.

\item  Notice that $E$ is of transcendence degree $1$ over ${\bf C}$ and, in
general, $E$ is not generated over ${\bf C}$ by a finite numbers of
elements. At the same time there exists a finite subset $A\subset E$ such
that every subfield $E^{\prime }$ of $E$ containing $A$ and invariant with
respect to $\limfunc{Aut}_{{\bf C}}E$ coincides with $E$. (Let $A$ be the
set of generators of ${\cal M}(X)$ over ${\bf C}$. Then $E$ is generated by $%
\bigcup\limits_{g\in \Gamma }gA$ over ${\bf C}$). Fields of this form were
studied in \cite{c5}.
\end{enumerate}
\end{description}

\fussy Let us consider the following situation. Let $G$ be a group, $U$ and $%
\Gamma $ its subgroups. Suppose that $G=\Gamma \cdot U$, $\Gamma \cap U=1$
(then $U\Gamma =(\Gamma U)^{-1}=G$). We get a bijection $\Gamma
\leftrightarrow G/U$. The group $G$ acts on $G/U$ by left translations, so $%
G $ acts on~$\Gamma $.\sloppy

\begin{lemma}
\label{le2.4}Let $A$ and $B$ be subsets of $\Gamma $.

\begin{enumerate}
\item  \label{le24.1}if $A$ is invariant with respect to $U$, then $A^{-1}$
is also invariant;

\item  \label{le24.2}$A$ is invariant\ $\Leftrightarrow UA\subset
AU\Leftrightarrow AU\subset UA\Leftrightarrow UA=AU$;

\item  \label{le24.3}if $A$ and $B$ are invariant, then $AB$ is also
invariant.
\end{enumerate}
\end{lemma}

\TeXButton{Proof}{\proof}Clearly $A$ is invariant $\Leftrightarrow UA\subset
AU$ and $A^{-1}$ is invariant $\Leftrightarrow UA^{-1}\subset
A^{-1}U\Leftrightarrow AU\subset UA$. So to prove \ref{le24.1}) and \ref
{le24.2}) it suffices to show that if $UA\subset AU$ than $AU\subset UA$.
Let $a\in A$, $\sigma \in U$. Then $a\sigma =\sigma ^{\prime }a^{\prime }$
for some $\sigma ^{\prime }\in U$, $a^{\prime }\in \Gamma $. Let us show
that $a^{\prime }\in A$. Indeed, $a^{\prime }=(\sigma ^{\prime
})^{-1}a\sigma \in UAU\subset AUU=AU$. So $a^{\prime }\in AU\cap \Gamma =A$.

To prove \ref{le24.3}) notice that if $UA\subset AU$ and $UB\subset BU$ then 
$UAB\subset AUB\subset ABU$.\TeXButton{End Proof}{\endproof}

\begin{lemma}
\label{le2.5}If $\Delta \subset \Gamma $ is a subgroup then the following
conditions are equivalent~:

\begin{enumerate}
\item  \label{le25.1}$U\cdot \Delta $ is a subgroup;

\item  \label{le25.2}$\Delta \cdot U$ is a subgroup;

\item  \label{le25.3}$\Delta $ is invariant with respect to the action of $U$%
.
\end{enumerate}
\end{lemma}

\TeXButton{Proof}{\proof}\ref{le25.1}) and \ref{le25.2}) are equivalent
because $\Delta U=(U\Delta )^{-1}$. \ref{le25.1}) means that $U\Delta
U\Delta \subset U\Delta $ and $(U\Delta )^{-1}\subset U\Delta $. Each of
these inclusions is equivalent to the inclusion $\Delta U\subset U\Delta $,
i.e., to \ref{le25.3}).\TeXButton{End Proof}{\endproof}

In Section \ref{Sec5} we shall need the following definition. Let 
\TeXButton{Discret}
{\discretionary{\mbox{$A_1,\ldots ,A_k,B_1,\ldots ,$}}{\mbox{$B_n\subset \Gamma $}}{\mbox{$A_1,\ldots ,A_k,B_1,\ldots ,B_n\subset \Gamma $}}}
be $U$-invariant subsets.

\begin{description}
\item[Definition.]  A relation of the form 
\begin{equation}
\label{e*}A_1\cdot \ldots \cdot A_k=B_1\cdot \ldots \cdot B_n
\end{equation}
is a $(k+n)$-tuple $(a_1,\ldots ,a_k,b_1,\ldots ,b_n)\in A_1\times \ldots
\times A_k\times B_1\times \ldots \times B_n$ such that $a_1\cdot \ldots
\cdot a_k=b_1\cdot \ldots \cdot b_n$.
\end{description}

Let us define an action of $U$ on the set of relations of the form (\ref{e*}%
). Let $\sigma \in U$ and $(a_1,\ldots ,a_k,b_1,\ldots ,b_n)$ be a relation
of the form (\ref{e*}). There exist unique $a_1^{\prime }\in A_1$, $\sigma
_1\in U$ such that $\sigma a_1=a_1^{\prime }\sigma _1$. There exist unique $%
a_2^{\prime }\in A_2$, $\sigma _2\in U$ such that $\sigma _1a_2=a_2^{\prime
}\sigma _2$, and so on. Thus we obtain $a_1^{\prime }\in A_1,\ldots
,a_k^{\prime }\in A_k$ and $\sigma _1,\ldots ,\sigma _k\in U$ such that $%
\sigma a_1\ldots a_k=a_1^{\prime }\ldots a_k^{\prime }\sigma _k$. Similarly,
we get $b_1^{\prime }\in B_1,\ldots ,b_n^{\prime }\in B_n$ and $\bar \sigma
_1,\ldots ,\bar \sigma _n\in U$ such that $\sigma b_1\ldots b_n=b_1^{\prime
}\ldots b_n^{\prime }\bar \sigma _n$. Since $a_1\ldots a_k=b_1\ldots b_n$,
we have $a_1^{\prime }\ldots a_k^{\prime }\sigma _k=b_1^{\prime }\ldots
b_n^{\prime }\bar \sigma _n$, Therefore $a_1^{\prime }\ldots a_k^{\prime
}=b_1^{\prime }\ldots b_n^{\prime }$. So $(a_1^{\prime },\ldots ,a_k^{\prime
},b_1^{\prime },\ldots ,b_n^{\prime })$ is a relation of the form (\ref{e*}%
). The action of $U$ on the set of relations of the form (\ref{e*}) is
defined as follows : $\sigma \in U$ maps $(a_1,\ldots ,a_k,b_1,\ldots ,b_n)$
to $(a_1^{\prime },\ldots ,a_k^{\prime },b_1^{\prime },\ldots ,b_n^{\prime
}) $. To show that this is really an action notice that for any $i\leq k$
and $j\leq n$ the product $a_1^{\prime }\ldots a_i^{\prime }$ is the result
of the action of $\sigma $ on $a_1\ldots a_i\in \Gamma =G/U$ and $%
b_1^{\prime }\ldots b_j^{\prime }$ is the result of the action of $\sigma $
on $b_1\ldots b_j$.

\section{The canonical diagram\label{Sec3}}

We shall need the following slight generalization of Artin's theorem in
Galois theory.

\begin{lemma}
\label{le3.1}Let $K$ be a Galois extension field of $F$, $U=Gal(K/F)$. Let $G
$ be a subgroup of $\limfunc{Aut}K$ such that $G\supset U$ and $[G:U]=n$.
Put $k=K^G$. Then $K$ is a Galois extension of $k$, $[F:k]=n$, $Gal(K/k)=G$.
\end{lemma}

This can be proved repeating word-for-word the arguments of Artin 
\cite[ch. VII, Theorem 2]{c6}.

\noindent {\bf Proof of Theorem \ref{th1.2}.}\qquad The diagram (\ref{e1.3})
with $\deg f=\limfunc{gcd}(\deg p,\deg q)$, $\deg g=(\deg p_1)\cdot (\deg
q_1)$ is unique if it exists, because ${\cal M}(W)={\cal M}(Y)\cap {\cal M}%
(Z)$ and ${\cal M}(V)$ is the composite of ${\cal M}(Y)$ and ${\cal M}(Z)$.
Indeed, \TeXButton{Discret}
{\discretionary{\mbox{$\deg p_1=\deg p/$}}{\mbox{$/\limfunc{gcd}(\deg p,\deg q)$}}{\mbox{$\deg p_1=\deg p/\limfunc{gcd}(\deg p,\deg q)$}}}
and $\deg q_1=\deg q/\limfunc{gcd}(\deg p,\deg q)$ are coprime, so ${\cal M}%
(V)$ is the composite of ${\cal M}(Y)$ and ${\cal M}(Z)$. Since $[{\cal M}%
(Y):{\cal M}(W)]=\deg g/\deg p_1=\deg q_1$ and $[{\cal M}(Z):{\cal M}%
(W)]=\deg p_1$ are coprime it follows that ${\cal M}(W)={\cal M}(Y)\cap 
{\cal M}(Z)$.

In the rest of the proof we use the notation of Section \ref{Sec2}. If $%
\Gamma $ is infinite, then $\Gamma $ is non-abelian (if $T_p$ and $T_q$
commute, then $\Gamma $ is finite). Obviously, in that case ${\cal M}(Y)\cap 
{\cal M}(Z)={\bf C}$.

If $\Gamma $ is finite, then $\Gamma _1$ is trivial because $J_1$ is
torsion-free. Therefore $\Gamma =\Gamma /\Gamma _1$ is cyclic of order $d=%
\limfunc{lcm}(\deg p,\deg q)$. By Lemma \ref{le2.3}, we have $[G:U]=d$. Put $%
F=E^G$. Using Lemma \ref{le3.1}, we obtain $[{\cal M}(X):F]=d$. Clearly, $%
{\cal M}(Y)\cap {\cal M}(Z)=F$. Let $W$ be the model of $F$. We get the
diagram : 
\begin{equation}
\label{e3.1}\TeXButton{Diagram 3.1 p.10 n.1}
{\matrix{\putCtrianglep<1`1`1;400>(-400,-400)[X`Y`W;p`r`p_0]
\putDtrianglep<0`1`1;400>(0,-400)[\phantom{X}`Z`\phantom{W};`q`q_0]}} 
\end{equation}

Put $\infty _W=r(\infty _X)$. Passing to the completions, we get :%
$$
\TeXButton{Diagram p.10 n.2}
{
\matrix{{\cal M}(X)_\infty&\mri{}&{\cal M}(Z)_\infty\cr
           \mdi{}&&\mdi{}\cr
           {\cal M}(Y)_\infty&\mri{}&F_\infty\cr}} 
$$

Since $\deg r=d$, it follows that $[{\cal M}(X)_\infty :F_\infty ]\leq d$.
Since $[{\cal M}(X)_\infty :{\cal M}(Y)_\infty ]$ and $[{\cal M}(X)_\infty :%
{\cal M}(Z)_\infty ]$ divide $[{\cal M}(X)_\infty :F_\infty ]$, it follows
that $d$ divides $[{\cal M}(X)_\infty :F_\infty ]$, hence $d=[{\cal M}%
(X)_\infty :F_\infty ]$. So $r:(X,\infty )\rightarrow (W,\infty )$ is a
generalized polynomial. Clearly, $p_0$ and $q_0$ are generalized polynomials
too. Besides, $T_p$ and $T_q$ generate a subgroup of $T_r$ of order $d$,
i.e., the group $T_r$.

Notice that $\deg p_0=\deg r/\deg p=\limfunc{lcm}(\deg p,\deg q)/\deg p$ and 
$\deg q_0=\limfunc{lcm}(\deg p,\deg q)/\deg q$ are coprime. So applying
Remark \ref{r0.1} to the diagram $(Y,\infty )\rightarrow (W,\infty
)\leftarrow (Z,\infty )$, one obtains the commutative diagram 
\begin{equation}
\label{e3.2}\TeXButton{Diagram 3.2 p.10 n.3}
{
\matrix{\xext=800 \yext=850
\adjust[(V,\infty)`;`;`;(W,\infty)`;]
\begin{picture}(\xext,\yext)(\xoff,\yoff)
\putCtrianglep<1`1`1;400>(0,0)[(V,\infty )`(Y,\infty )`(W,\infty );p_1`g`p_0]
\putDtrianglep<0`1`1;400>(400,0)[\phantom{X}`(Z,\infty )`\phantom{W};`q_1`q_0]
\end{picture}
}} 
\end{equation}
where $V$ is the normalization of $Y\times _WZ$, $\deg q_1=\deg p_0$, $\deg
p_1=\deg q_0$. From the definition of $Y\times _WZ$ it follows that the
diagrams (\ref{e3.1}) and (\ref{e3.2}) can be included into a diagram of the
form (\ref{e1.3}). Clearly, this is the desired diagram.\TeXButton{End Proof}
{\endproof}

\begin{description}
\item[Remark.]  We have two equivalence relation on $X$ : $%
R_p=\{(x_1,x_2)\in X\times X\mid p(x_1)=p(x_2)\}$, $R_q=\{(x_1,x_2)\in
X\times X\mid q(x_1)=q(x_2)\}$. Denote by $R$ the equivalence relation
generated by them. If $\Gamma $ is finite it is easy to show that $R$ is an
algebraic curve on $X\times X$. The essential part of the proof of Theorem 
\ref{th1.2} is the construction of the quotient $X/R$ as a Riemann surface.
We have done it using Lemma \ref{le3.1} (in fact, we have constructed the
field ${\cal M}(X/R)$). One can also construct the Riemann surface $X/R$
directly using Theorem G from \cite[Appendix A]{c7}. Besides, one can
construct the algebraic curve $X/R$ in the framework of algebraic geometry
using Theorem 4.1 from \cite[expos\'e V, p.262]{c8}.

\item[Corollary.]  Let $p_i:(X,\infty )\rightarrow (Y_i,\infty )$ be
generalized polynomials, $i\in \{1,2,3\}$. We have the following diagram of
fields :%
$$
\TeXButton{Diagram p.11 n.1}
{
\matrix{{\cal M}(Y_1)&\mli{}&{\cal M}(X)&\mri{}&{\cal M}(Y_3)\cr
           &&\mdi{}&&\cr
           &&{\cal M}(Y_2)&&\cr}} 
$$
\fussy Suppose ${\cal M}(Y_i)\cap {\cal M}(Y_j)\neq {\bf C}$ for any $i,j$;
then ${\cal M}(Y_1)\cap {\cal M}(Y_2)\cap {\cal M}(Y_3)\neq {\bf C}$.%
\sloppy 
\end{description}

\TeXButton{Proof}{\proof}By Theorem \ref{th1.2} the elements $T_{p_i}$ and $%
T_{p_j}$ commute, and we have a diagram of generalized polynomials%
$$
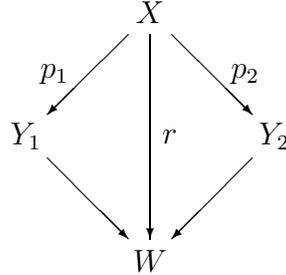
\TeXButton{Diagram p.11 n.2}
{
\matrix{\xext=800 \yext=850
\adjust[X`;`;`;W`;]
\begin{picture}(\xext,\yext)(\xoff,\yoff)
\putCtrianglep<1`1`1;400>(0,0)[X`Y_1`W;p_1`r`]
\putDtrianglep<0`1`1;400>(400,0)[\phantom{X}`Y_2`\phantom{W};`p_2`]
\end{picture}
}} 
$$
such that ${\cal M}(W)={\cal M}(Y_1)\cap {\cal M}(Y_2)$, $T_{p_1}$ and $%
T_{p_2}$ generate $T_r$. The elements of $T_r$ and $T_{p_3}$ commute.
Applying Theorem \ref{th1.2} to the pair $r,p_3$, we get ${\cal M}(Y_1)\cap 
{\cal M}(Y_2)\cap {\cal M}(Y_3)\neq {\bf C}$.\TeXButton{End Proof}{\endproof}

Proposition \ref{pp1.1} follows immediately from this corollary. So this
proposition is also proved.

Suppose we are given a commutative diagram of generalized polynomials :%
$$
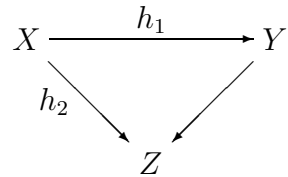
\Vtrianglep<1`1`1;400>[X`Y`Z;h_1`h_2`] 
$$
Then $T_{h_1}$ is a subgroup of $T_{h_2}$. A right factor $h_1$ of $h_2$ can
be reconstructed from the group $T_{h_1}$ as follows : ${\cal M}(Y)={\cal M}%
(X)\cap {\cal M}(X)_\infty ^{T_{h_1}}$. Thus certain subgroups of $T_{h_2}$
corresponds to right factors of $h_2$ (not necessarily all the subgroups!).

Consider a diagram (\ref{e1.1}) again. In Section \ref{Sec1} the
intermediate field ${\cal M}(Z)\subset F_{q,p}\subset {\cal M}(X)$ was
introduced. The following theorem indicates the subgroup of $T_q$
corresponding to this field. There exists a commutative diagram of
generalized polynomials :%
$$
\TeXButton{Diagram p.12 n.2}
{
    \matrix{X&\mr{q_1}&X_{q,p}&\mr{q_2}&Z\cr
               \md{p}&&\md{}&&&\cr
               Y&\mr{}&Y_{q,p}&&&\cr}
} 
$$
such that $q_2\circ q_1=q$, ${\cal M}(X_{q,p})=F_{q,p}$, ${\cal M}%
(Y_{q,p})=F_{q,p}\cap {\cal M}(Y)$. Put $H_{q,p}=\left\{ \sigma \in T_q\mid
\tau \sigma =\sigma \tau \text{ for all }\tau \in T_p\right\} $.

\begin{theorem}
\label{th3.1}$T_{q_1}=H_{q,p}$.
\end{theorem}

To prove the theorem we need several lemmas.

\begin{lemma}
\label{le3.2}Let $H$ be a finite subgroup of $J$. Then $H$ is a cyclic group.
\end{lemma}

\TeXButton{Proof}{\proof}Since $J_1$ is torsion-free, it follows that $H\cap
J_1=1$. So $H\simeq H/H_1\hookrightarrow J/J_1\simeq {\bf C}^{*}$. A finite
subgroup of ${\bf C}^{*}$ is cyclic.\TeXButton{End Proof}{\endproof}

We use the constructions and notation of Section \ref{Sec2}.

\begin{lemma}
\label{le3.3}Let $L\subset \Gamma $ be a subset invariant under the action
of $U$, i.e., $uLU=LU$ for all $u\in U$. Put 
$$
N_L=\left\{ \sigma \in T_q\mid L\sigma =L\right\} . 
$$
Then $N_L$ is a subgroup of $T_q$ invariant under the action of $U$.
\end{lemma}

\begin{description}
\item[Remark.]  The subgroups $T_p$ and $T_q$ of $\Gamma $ are invariant
under the action of $U$.
\end{description}

\TeXButton{Proof}{\proof}It is easily checked that $N_L$ is a subgroup.
According to Lemma \ref{le2.4}, it remains to show that $UN_L\subset
N_L\cdot U$. Suppose $\sigma \in N_L$, $\gamma \in U$, $\gamma \sigma
=\sigma ^{\prime }\gamma ^{\prime }$, where $\sigma ^{\prime }\in T_q$, $%
\gamma ^{\prime }\in U$. We must prove that $L\sigma ^{\prime }=L$. Using
Lemma \ref{le2.4}, we have $LU\subset UL$, hence $LUN_L\subset ULN_L\subset
UL\subset LU$. Therefore $L\sigma ^{\prime }=L\gamma \sigma (\gamma ^{\prime
})^{-1}\subset LU$. On the other hand, $L\sigma ^{\prime }\subset \Gamma $
and $\Gamma \cap LU=L$, hence $L\sigma ^{\prime }\subset L$. Since $\sigma
^{\prime }$ is of finite order, it follows that $L\sigma ^{\prime }=L$.%
\TeXButton{End Proof}{\endproof}

\noindent {\bf Proof of Theorem \ref{th3.1}.\qquad }Put $L=T_qT_p$. By \ref
{le24.3}) of Lemma \ref{le2.4}, $L$ is invariant under the action of $U$. By
Lemma \ref{le3.3}, we get a subgroup $N_L(=N)$ of $T_q$ such that $N$ is
invariant with respect to $U$. By Lemma \ref{le2.5}, $N\cdot U$ is a
subgroup of $G$, $[NU:U]=\#N$. Further, we have ${\cal M}(Z)\subset
E^{NU}\subset {\cal M}(X)$, $[{\cal M}(X):E^{NU}]=\#N$. This means that $N$
corresponds to a right factor of $q$. Obviously, $H_{q,p}\subset N$. Let us
show that $H_{q,p}=N$.

We have $T_pN\subset T_qT_p$. Let $\tau \sigma =\sigma ^{\prime }\tau
^{\prime }$, where $\tau ,\tau ^{\prime }\in T_p$, $\sigma \in N$, $\sigma
^{\prime }\in T_q$. Then $T_qT_p\sigma ^{\prime }=T_qT_p\sigma (\tau
^{\prime })^{-1}=T_qT_p(\tau ^{\prime })^{-1}=T_qT_p$. So $\sigma ^{\prime
}\in N$. We get the property $T_pN\subset NT_p$. This implies that $N\cdot
T_p$ is a subgoup of $\Gamma $. By Lemma \ref{le3.2}, $N\cdot T_p$ is
abelian. Therefore $N\subset H_{q,p}$.

We have proved that $H_{q,p}$ corresponds to a right factor of $q$. Using
Theorem \ref{th1.2}, it is not hard to check that $H_{q,p}$ corresponds to $%
F_{q,p}$.\TeXButton{End Proof}{\endproof}

By definition, the pair $p,q$ in (\ref{e1.1}) is irreducible if $%
F_{q,p}=F_{p,q}={\cal M}(X)$. By Theorem \ref{th3.1}, this is equivalent to
the property $H_{p,q}=H_{q,p}=1$.

\begin{description}
\item[Example.]  Consider the following pair of polynomials : $p(z)=z^n$, $%
q(z)=(z+1)^m$. We have $T_p=\left\{ z\mapsto \varepsilon z\mid \varepsilon
^n=1\right\} $, $T_q=\left\{ z\mapsto \delta z+(\delta -1)\mid \delta
^m=1\right\} $. Now it is easy to check that $H_{p,q}=H_{q,p}=1$. So this
pair is irreducible.
\end{description}

\begin{lemma}
\label{le3.4}Consider a diagram (\ref{e1.1}) such that ${\cal M}(Y)\cap 
{\cal M}(Z)\neq {\bf C}$. Let $K$ be an intermediate field : $({\cal M}%
(Y)\cap {\cal M}(Z))\subset K\subset {\cal M}(Y)$. Let $\tilde K$ be the
composite of $K$ and ${\cal M}(Z)$. Then $\tilde K\cap {\cal M}(Y)=K$.
\end{lemma}

\TeXButton{Proof}{\proof}Let $F={\cal M}(Y)\cap {\cal M}(Z)$. We get the
following commutative diagram%
$$
\TeXButton{Diagram p.14 n.1}
{
\matrix{{\cal M}(X)&\mri{}&{\tilde K}&\mri{}&{\cal M}(Z)\cr
            \mdi{}&&\mdi{}&&\mdi{}\cr
            {\cal M}(Y)&\mri{}&K&\mri{}&F\cr}} 
$$
By Theorem \ref{th1.2}, $[{\cal M}(Y):F]$ and $[{\cal M}(Z):F]$ are coprime.
Therefore $[K:F]$ and $[{\cal M}(Z):F]$ are coprime, hence $[\tilde K:K]=[%
{\cal M}(Z):F]$. It follows that $[{\cal M}(Y):K]$ and $[\tilde K:K]$ are
coprime, so $K=\tilde K\cap {\cal M}(Y)$.\TeXButton{End Proof}{\endproof}

\noindent {\bf Proof of Proposition \ref{pp1.2}.}\quad \quad It follows
immediately from the previous lemma.\TeXButton{End Proof}{\endproof}

\noindent {\bf Proof of Proposition \ref{pp1.3}.}\quad \quad The assertion
is nontrivial for the formal discreteness property. Taking into account
Theorem \ref{th1.2}, it suffices to prove the following result.

\begin{lemma}
\label{le3.5}Suppose we are given the following diagram of generalized
polynomials :%
$$
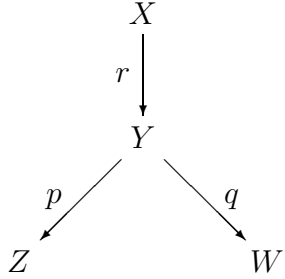
\TeXButton{Diagram p.14 n.2}
{\matrix{\xext=800 \yext=850
\adjust[X`;`;`;W`;]
\begin{picture}(\xext,\yext)(\xoff,\yoff)
\putmorphism(400,800)(0,-1)[X``r]{400}1l
\putAtriangle<1`1`0;400>(0,0)[Y`Z`W;p`q`]
\end{picture}
}
} 
$$
Let $\Gamma $ be the group generated by $T_p$ and $T_q$, $\Gamma ^{\prime }$
the group generated by $T_{p\circ r}$ and $T_{q\circ r}$. Then $\Gamma $ is
formally discrete iff $\Gamma ^{\prime }$ is formally discrete.
\end{lemma}

\TeXButton{Proof}{\proof}Put $J(X,Y,\infty )=\{(g_X,g_Y)\mid g_X\in
J(X,\infty ),g_Y\in J(Y,\infty ),g_Y\circ r=r\circ g_X\}$. $J(X,Y,\infty )$
is a subgroup of $J(X,\infty )\times J(Y,\infty )$. Consider the projections 
$\pi :J(X,Y,\infty )\rightarrow J(Y,\infty )$, $j:J(X,Y,\infty )\rightarrow
J(X,\infty )$. It is easy to see that $\pi $ is surjective, $\limfunc{Ker}%
\pi =T_r$, $j$ is injective, $T_{p\circ r}=j(\pi ^{-1}(T_p))$, $T_{q\circ
r}=j(\pi ^{-1}(T_q))$. So $\Gamma ^{\prime }=j(\pi ^{-1}(\Gamma ))$.

Choose a meromorphic function $z$ in a neighborhood of $\infty \in Y$ with a
pole of order $1$ at $\infty $ (so $z^{-1}$ is a local coordinate at $\infty 
$). Choose a similar function $\zeta $ in a neighborhood of $\infty \in X$
so that $r^{*}(z)=\zeta ^n$, $n=\deg r$ (i.e., in terms of $z$ and $\zeta $
the mapping $r$ is described by $z=\zeta ^n$). Then we can write $g_X\in
J(X,\infty )$ and $g_Y\in J(Y,\infty )$ as $g_X(\zeta
)=\sum\limits_{j=-1}^\infty a_j\zeta ^{-j}$, $g_Y(z)=\sum\limits_{k=-1}^%
\infty b_kz^{-k}$, $a_{-1}\neq 0$, $b_{-1}\neq 0$. Further, the relation $%
g_Y\circ r=r\circ g_X$ can be written as $\left( \sum\limits_{j=-1}^\infty
a_j\zeta ^{-j}\right) ^n=\sum\limits_{k=-1}^\infty b_k\zeta ^{-kn}$ or $%
\sum\limits_{j=-1}^\infty a_j\zeta ^{-j}=\zeta (b_{-1}+b_0\zeta
^{-n}+b_1\zeta ^{-2n}+\ldots )^{\frac 1n}$. Now it is clear that a subgroup $%
\Gamma \subset J(Y,\infty )$ is formally discrete iff $j(\pi ^{-1}(\Gamma ))$
is formally discrete.\TeXButton{End Proof}{\endproof}

Proposition \ref{pp1.3} is also proved.\TeXButton{End Proof}{\endproof}

\section{Irreducible pairs of generalized polynomials\label{Sec4}}

\begin{theorem}
\label{th4.1}Suppose we are given a diagram (\ref{e1.1}) such that the pair $%
p,q$ is irreducible and $\deg p>1$, $\deg q>1$. Let $\Gamma $ be the group
generated by $T_p$ and $T_q$. Suppose $\Gamma _1$ is abelian. Put $n=\deg p$%
, $m=\deg q$. Put $\tilde p(z)=z^n$, $\tilde q(z)=(z+1)^m$, $\tilde p,\tilde
q\in {\bf C}[z]$. Let $\tilde \Gamma $ be the group generated by $T_{\tilde
p}$ and $T_{\tilde q}$. Then there exists an isomorphism $\varphi :\Gamma
\tilde \rightarrow \tilde \Gamma $ such that $\varphi |_{T_p}:T_p\tilde
\rightarrow T_{\tilde p}$, $\varphi |_{Tq}:T_q\tilde \rightarrow T_{\tilde q}
$, and $\varphi |_{\Gamma _1}:\Gamma _1\tilde \rightarrow \tilde \Gamma _1$
are isomorphisms. Besides, $\Gamma $ is formally discrete iff $\tilde \Gamma 
$ is formally discrete.
\end{theorem}

To prove the theorem we need the following result.

Let $k\geq 1$ be an integer. Denote by $g_{z^{k+1}}^t$ the germ $({\bf C}%
,0)\rightarrow ({\bf C},0)$ of a time $t$ map for the flow of the
holomorphic vector field $z^{k+1}\frac d{dz}$. The set of germs $%
G(k)=\{\lambda g_{z^{k+1}}^t\mid \lambda \in {\bf C}^{*},t\in {\bf C\}}$ is
a group with respect to superposition. For brevity, denote $\lambda
g_{z^{k+1}}^t$ by $(\lambda ,t)$. The multiplication table for $G(k)$ has
the following form : 
\begin{equation}
\label{e**}(\lambda ,t)\times (\mu ,s)=(\lambda \mu ,t\mu ^k+s). 
\end{equation}
The subgroup $C(k)=\{\lambda \in {\bf C}\mid \lambda ^k=1\}$ is the center
of $G(k)$. Put $G_d(k)=\{\lambda g_{z^{k+1}}^t\in G(k)\mid \lambda ^d=1\}$.
Then $G_d(k)$ is a subgroup of $G(k)$. It is easy to see that if $h\in
G_k(k) $ is an element of finite order, then $h\in C(k)$.

\begin{description}
\item[Theorem A.]  {\bf (\cite[Theorem 2.2, p.66]{c9})} A finitely generated
nonabelian solvable group of germs of conformal mappings $({\bf C}%
,0)\rightarrow ({\bf C},0)$ is formally equivalent to a finitely generated
subgroup of the group $G(k)$ for some $k$.
\end{description}

\noindent {\bf Proof of Theorem \ref{th4.1}.}\quad \quad Put $d=\limfunc{lcm}%
(n,m)$. Choose a local parameter $z$ at $\infty \in X$ and identify the
group $J(X,\infty )$ with the group of germs of conformal mappings : $({\bf C%
},0)\rightarrow ({\bf C},0)$. By Remark \ref{r1.4}, ${\cal M}(Y)\cap {\cal M}%
(Z)={\bf C}$. By Theorem \ref{th1.2}, $\Gamma $ is nonabelian. One the other
hand, $\Gamma _1$ is abelian, hence $\Gamma $ is solvable. By Theorem A, $%
\Gamma $ is formally equivalent to a subgroup of $G_d(k)$ for some $k$. We
have $T_p\cap G_k(k)\subset C(k)$, $T_q\cap G_k(k)\subset C(k)$. By Theorem 
\ref{th3.1}, $H_{p,q}=H_{q,p}=1$. It follows that $T_p\cap G_k(k)=T_q\cap
G_k(k)=1$. Therefore $\limfunc{gcd}(n,k)=\limfunc{gcd}(m,k)=1$, so $\limfunc{%
gcd}(d,k)=1$. Let the map $f:G_d(k)\rightarrow G_d(1)$ be given by $\lambda
g_{z^{k+1}}^t\mapsto \lambda ^kg_{z^2}^t$. Since $\gcd (d,k)=1$, it follows
that $f$ is bijective. The multiplication table (\ref{e**}) shows that $f$
is an isomorphism.

Let $h_p$ (resp. $h_q$) be a generator of $T_p$ (resp. $T_q$). Let $%
\varepsilon g_{z^2}^{t_1}$, $\delta g_{z^2}^{t_2}$ be the elements of $%
G_d(1) $ corresponding to $h_p$ and $h_q$ respectively (Then $\varepsilon $
(resp. $\delta $) is a primitive $n$-th (resp. $m$-th) root of unity.).
Using the conjugation in $G_d(1)$, we may assume that $t_1=0$. Then $t_2\neq
0$. Notice that the map given by $\lambda g_{z^2}^t\mapsto \lambda
g_{z^2}^{ct}$ ($c\in {\bf C}^{*}$) is an automorphism of $G_d(1)$. Finally,
the pair of generators becomes $\varepsilon ,\delta g_{z^2}^1$. Clearly, $%
\Gamma ^{\prime }$ is formally discrete iff $\Gamma $ is formally discrete.%
\TeXButton{End Proof}{\endproof}

\begin{lemma}
\label{le4.1}Let $p(z)=z^n$, $q(z)=(z+1)^m$; $n,m\geq 2$. Let $\Gamma $ be
the group generated by $T_p$ and $T_q$. The group $\Gamma $ is formally
discrete iff $\limfunc{lcm}(n,m)\in \{2,3,4,6\}$.
\end{lemma}

\TeXButton{Proof}{\proof}Let $\varepsilon $ (resp. $\delta $) be a primitive 
$n$-th (resp. $m$-th) root of unity. We have $\varepsilon \Gamma _1=\Gamma
_1 $, $\delta \Gamma _1=\Gamma _1$. (Here $\Gamma _1$ is considered as a
subgroup of ${\bf C}$). Therefore $e^{2\pi i/d}\Gamma _1\subset \Gamma _1$,
where $d:=\limfunc{lcm}(n,m)$. So if $\Gamma _1$ is discrete then $d\in
\{2,3,4,6\}$. On the other hand, $\Gamma _1\subset {\bf Z}[e^{2\pi i/d}]$,
so if $d\in \{2,3,4,6\}$ then $\Gamma _1$ is discrete and $\Gamma $ is
formally discrete.\TeXButton{End Proof}{\endproof}

\begin{description}
\item[Corollary.]  Under the conditions of Theorem \ref{th4.1} the group $%
\Gamma $ is formally discrete iff $\limfunc{lcm}(n,m)\in \{2,3,4,6\}$.
\end{description}

\begin{lemma}[Main group-theoretic lemma]
\label{MGTL}Put $p(z)=z^n$, $q(z)=(z+1)^m$, $n,m\geq 2$. Let $\Gamma $ be
the group generated by $T_p$ and $T_q$. Let $G$ be an abstract group and $U$
its subgroup. Suppose $\Gamma $ is embedded into $G$ as a subgroup and $%
\Gamma U=G$, $\Gamma \cap U=1$. Suppose $T_pU$ and $T_qU$ are subgroups of $G
$. Suppose that $(n,m)\in P_1\cup P_2\cup P_3$, where $P_1=\{(n,m)\mid n=m\}$%
, $P_2=\{(n,m)\mid n=2$ or $m=2\}$, and $P_3$ consists of $(3,6)$ and $(6,3)$%
. Then there exists a subgroup $U^{\prime }$ of $U$ such that $[U:U^{\prime
}]<\infty $ and $U^{\prime }$ is a normal subgroup of $G$.
\end{lemma}

This lemma will be proved in the following Section.

\begin{description}
\item[Remark.]  If $\limfunc{lcm}(n,m)\in \{2,3,4,6\}$ and $n,m\geq 2$ then $%
(n,m)\in P_1\cup P_2\cup P_3$.
\end{description}

\begin{theorem}
\label{th4.2}Suppose we are given a diagram (\ref{e1.1}) such that the pair $%
p,q$ is irreducible and $\deg p>1$, $\deg q>1$. Let $\Gamma $ be the group
generated by $T_p$ and $T_q$. Suppose $\Gamma _1$ is abelian. Put $n=\deg p$%
, $m=\deg q$, $\tilde p(z)=z^n$, $\tilde q(z)=(z+1)^m$. If $(n,m)$ belongs
to the set $P_1\cup P_2\cup P_3$ from the main group-theoretic lemma then
there exists a commutative diagram :%
$$
\TeXButton{Diagram p.18 n.1}
{
\matrix{(Y,\infty)&\ml{p}&(X,\infty)&\mr{q}&(Z,\infty)\cr
           \md{}&&\md{}&&\md{}\cr
           (\bf{CP}^1,\infty)&\ml{\tilde p}&(\bf{CP}^1,\infty)&\mr{\tilde q}&(\bf{CP}^1,\infty)\cr}}
$$
where the vertical arrows are isomorphisms.
\end{theorem}

\TeXButton{Proof}{\proof}We use the notation of Section \ref{Sec2}. By
Theorem \ref{th4.1}, we can apply the main group-theoretic lemma. We get a
subgroup $U^{\prime }$ of $U$ such that $[U:U^{\prime }]<\infty $ and $%
U^{\prime }$ is a normal subgroup of $G$. The field $E^{U^{\prime }}$ is
normal over both ${\cal M}(Y)$ and ${\cal M}(Z)$, hence $E^{U^{\prime }}=E$,
i.e., $U^{\prime }=1$. So $\#U<\infty $ and we get a diagram : 
$$
\TeXButton{Diagram p.18 n.2}
{\matrix{\xext=800 \yext=850
\adjust[X`;`;`;W`;]
\begin{picture}(\xext,\yext)(\xoff,\yoff)
\putmorphism(400,800)(0,-1)[W``r]{400}1l
\putAtriangle<1`1`0;400>(0,0)[X`Y`Z;p`q`]
\end{picture}
}
} 
$$
where $W$ is a compact Riemann surface, ${\cal M}(W)=E$, $r$ is nonconstant
holomorphic, $p\circ r$ and $q\circ r$ are Galois coverings. Since $G\subset 
\limfunc{Aut}W$, we have $\#\limfunc{Aut}W=\infty $. Therefore, $W$ is of
genus $0$ or $1$. $G$ acts on the finite set $S=r^{-1}(\infty )\subset W$.
If $W$ is of genus $1$, then for any $w\in W$ the group $\{g\in \limfunc{Aut}%
W\mid gw=w\}$ is finite. Therefore the group $\{g\in \limfunc{Aut}W\mid
gS=S\}$ is finite too. So $W$ is of genus $0$. Put $G_0=\{g\in G\mid \forall
s\in S:g(s)=s\}$. Then $[G:G_0]<\infty $, hence $\#G_0=\infty $. It means
that $\#S\leq 2$. By Lemma \ref{le2.2}, $r$ is not ramified over $\infty \in
X$. Suppose $\#S=2$. By \ref{le23.2}) of Lemma \ref{le2.3}, $\Gamma =G_0$.
It can be assumed that $W={\bf CP}^1$, $S=\{0,\infty \}$. We have $\Gamma
\subset \{g\in \limfunc{Aut}{\bf CP}^1\mid g(0)=0,g(\infty )=\infty \}\simeq 
{\bf C}^{*}$, hence $\Gamma $ is abelian, which is impossible. Thus $\#S=1$,
i.e., $r$ is an isomorphism, $p$ and $q$ are Galois coverings. This
completes the proof.\TeXButton{End Proof}{\endproof}

Theorem \ref{th1.3} is a special case of Theorem \ref{th4.2}.

\section{Proof of the main group-theoretic lemma\label{Sec5}}

The main group-theoretic lemma was formulated in Section \ref{Sec4}.

\begin{description}
\item[Remark.]  The basic idea of our proof is to study the action of $U$ on
the set of relations of some form in $\Gamma $ (this action was defined at
the end of Section \ref{Sec2}). Put $A=T_p\setminus \{id\}$, $B=T_q\setminus
\{id\}$. We use the relations of the form 
\begin{equation}
\label{e5.1}B\cdot A=A\cdot B
\end{equation}
and of the form 
\begin{equation}
\label{e5.2}A\cdot B\cdot A=B\cdot A\cdot B.
\end{equation}
\end{description}

\begin{proposition}
\label{pp5.1}The main group-theoretic lemma holds for $n=m$.
\end{proposition}

\TeXButton{Proof}{\proof}Put $G_p=T_pU$, $G_q=T_qU$, $U_p=\left\{ \sigma \in
U\mid \sigma \tau U=\tau U\text{ for all }\tau \in T_p\right\} $, $%
U_q=\left\{ \sigma \in U\mid \sigma \tau U=\tau U\text{ for all }\tau \in
T_q\right\} $. Clearly, $U_p$ is the kernel of the left action of $G_p$ on $%
G_p/U$. So $U_p$ is a normal subgroup of $G_p$ and $[U:U_p]<\infty $. Since $%
U_q$ is a normal subgroup of $G_q$, it suffices to show that $U_p=U_q$ (this
will imply that $U_p$ is a normal subgroup of $G$).

Put $h_p(z)=\varepsilon z$, $h_q(z)=\varepsilon z+(\varepsilon -1)$. Here $%
h_p\in T_p$, $h_q\in T_q$, $\varepsilon $ is a primitive $n$-th root of
unity. Let us consider two cases.

\begin{enumerate}
\item[\bf{Case 1}]  $n$ is even.

\begin{enumerate}
\item  If $n=2$, then $U=U_p=U_q$, there is nothing to prove.

\item  Assume $n\geq 4$. Consider the set of relations of the form (\ref
{e5.1}). If $h_q^{l_1}h_p^{l_2}=h_p^{k_1}h_q^{k_2}$, then $\varepsilon
^{l_1+l_2}=\varepsilon ^{k_1+k_2}$ and $\varepsilon ^{l_1}-1=\varepsilon
^{k_1}(\varepsilon ^{k_2}-1)$, so it is easily checked that the relation of
the form (\ref{e5.1}) are precisely the following ones : 
\begin{equation}
\label{e5.3}h_q^lh_p^{-l+\frac n2}=h_p^{l+\frac n2}h_q^{-l}\text{,}
\end{equation}
where $2l\not \equiv 0{\rm \;mod\;}n$. Notice that all $h_p^s$ ($s\not
\equiv 0{\rm \;mod\;}\frac n2$) occur in the right hand side of (\ref{e5.3}%
). Let us show that $U_p=U_q$. Let $\sigma \in U_q$. Then $\sigma $
preserves each relation (\ref{e5.3}), hence $\sigma $ preserves $h_p^i$ for
every $i\not \equiv 0{\rm \;mod\;}\left( \frac n2\right) $. Since $\sigma $
preserves $T_p$ and $id$, $\sigma $ preserves $h_p^{\frac n2}$. So $\sigma
\in U_p$. Similarly, $U_p\subset U_q$.
\end{enumerate}

\item[\bf{Case 2}]  $n$ is odd.\\Since $\Gamma _1$ is abelian, $h_q^{-l}h_p^l
$ and $h_q^{-s}h_p^s$ commute. This provides the following relations : 
\begin{equation}
\label{e5.4}h_p^{l_1}h_q^{l_2}h_p^{l_3}=h_q^{-l_3}h_p^{-l_2}h_p^{-l_1}\text{,%
}
\end{equation}
where $l_i\not \equiv 0{\rm \;mod\;}n$, $l_1+l_2+l_3\equiv 0{\rm {\;mod\;}}n$%
. The relations of the form (\ref{e5.2}) contain the relations (\ref{e5.4}).
We shall not find all the relations of the form (\ref{e5.2}), but prove the
following lemma.
\end{enumerate}

\begin{lemma}
\label{le5.1}Suppose that 
\begin{equation}
\label{e5.5}h_p^{l_1}h_q^{l_2}h_p^{l_3}=h_p^{k_1}h_q^{k_2}h_p^{k_3},\qquad
l_i,k_i\TeXButton{Not}{\not\equiv}0{\rm {\;mod\;}}n{\rm \text{.}}
\end{equation}
Then $l_1\not \equiv k_1{\rm {\;mod\;}}n$ (recall that $n$ is odd!).
\end{lemma}

\noindent {\bf Deduction of Case 2 from Lemma \ref{le5.1}.}\qquad From (\ref
{e5.4}) it follows that for any $l_1\not \equiv 0{\rm {\;mod\;}}n$, $k_1\not
\equiv 0{\rm \;mod\;}n$, $k_1\not \equiv l_1{\rm \;mod\;}n$ there exists a
relation of the form (\ref{e5.5}) with these $l_1,k_1$. Let $\sigma \in U$, $%
\sigma h_p^sU=h_p^tU$ for some $s,t$. Then $\sigma $ takes the set of
relations (\ref{e5.5}) such that $h_p^{l_1}=h_p^s$ to the set of those
relations that $h_p^{l_1}=h_p^t$. Therefore $\sigma $ takes the set $\left\{
h_q^j\mid j\not \equiv s{\rm \;mod\;}n\right\} $ to the set $\left\{
h_q^j\mid j\not \equiv t{\rm \;mod\;}n\right\} $, hence $\sigma
h_q^sU=h_q^tU $. Finally, for $\sigma \in U$ we have $\sigma h_q^sU=h_q^tU$
iff $\sigma h_p^sU=h_p^tU$. This implies $U_p=U_q$.\TeXButton{End Proof}
{\endproof}

It remains to prove Lemma \ref{le5.1}.

Put $\alpha _l=h_q^lh_p^{-l}$. Then $\alpha _l(z)=z+(\varepsilon ^l-1)$.
Suppose $h_p^{l_1}h_q^{l_2}h_p^{l_3}=h_q^{l_1}h_p^{k_2}h_q^{k_3}$, $l_i,k_i%
\not \equiv 0{\rm \;mod\;}n$. Then $k_2+k_3=l_2+l_3$. Therefore $%
(h_q^{-l_1}h_p^{l_1})\cdot (h_q^{l_2}h_p^{-l_2})=(h_p^{k_2}h_q^{-k_2})\cdot
(h_q^{k_2+k_3}h_p^{-k_2-k_3})$, i.e., $(\varepsilon ^{-l_1}-1)+(\varepsilon
^{l_2}-1)=(\varepsilon ^{k_2+k_3}-1)-(\varepsilon ^{k_2}-1)$. So it suffices
to prove the following lemma.

\begin{lemma}
\label{le5.2}Let $\varepsilon $ be a primitive $n$-th root of unity. Let $n$
be odd. Then the equation 
\begin{equation}
\label{e5.6}\varepsilon ^{\gamma _1}+\varepsilon ^{\gamma _2}+\varepsilon
^{\gamma _3}=\varepsilon ^\mu +2 
\end{equation}
has no solutions such that $\gamma _1,\gamma _2,\gamma _3\not \equiv 0{\rm %
\;mod\;}n$ (here $\gamma _i,\mu _i\in {\bf Z}$ are unknowns).
\end{lemma}

\TeXButton{Proof}{\proof} If $\mu \equiv 0{\rm \;mod\;}n$, then $\varepsilon
^{\gamma _1}=\varepsilon ^{\gamma _2}=\varepsilon ^{\gamma _3}=1$, i.e., a
contradiction. Therefore $\mu \not \equiv 0{\rm \;mod\;}n$.

The idea is to average by the action of $Gal({\bf Q}(\varepsilon )/{\bf Q})$.

Put $K=\bigcup\limits_m{\bf Q}(\sqrt[m]{1})$. For any $m\in {\bf N}$ define
a ${\bf Q}$-linear functional $T_m:{\bf Q}(\sqrt[m]{1})\rightarrow {\bf Q}$
by $T_m(z)=\frac 1{\#H}\sum\limits_{h\in H}h(z)$, where $H=Gal({\bf Q}(\sqrt[%
m]{1})/{\bf Q})$. If $m^{\prime }\mid m$, then $T_m|_{{\bf Q}(\sqrt[%
m^{\prime }]{1})}=T_{m^{\prime }}$. Therefore a ${\bf Q}$-linear functional $%
T:K\rightarrow {\bf Q}$ is well defined by $T|_{{\bf Q}(\sqrt[m]{1})}=T_m$.
Recall that if $\delta $ is a primitive $m$-th root of unity, then $Tr_{{\bf %
Q}}^{{\bf Q}(\delta )}(\delta )=\mu (m)$, where $\mu $ is the M\"obius
function. Therefore, $T\delta =\frac{\mu (m)}{\varphi (m)}$. One has $%
\varphi (1)=1=\varphi (2)$, $\varphi (3)=\varphi (4)=\varphi (6)=2$, $%
\varphi (m)>2$ for $m\neq 1,2,3,4,6$. So if $m>1$ is odd then $-\frac 12\leq
T\delta <\frac 12$.

Applying $T$ to (\ref{e5.6}), we obtain%
$$
\frac 32\leq 2+T(\varepsilon ^\mu )=T(\varepsilon ^{\gamma
_1})+T(\varepsilon ^{\gamma _2})+T(\varepsilon ^{\gamma _3})<\frac 32\text{,}
$$
i.e., a contradiction.\TeXButton{End Proof}{\endproof}

Proposition \ref{pp5.1} is proved.

\begin{proposition}
\label{pp5.2}The main group-theoretic lemma holds for $n=3$, $m=6$.
\end{proposition}

\TeXButton{Proof}{\proof}Let $\omega $ be a primitive $6$-th root of unity.
Put $h_p(z)=\omega ^2z$, $h_q(z)=\omega z+(\omega -1)$. It is easy to show
that there are exactly two relations of the form~(\ref{e5.1}) 
\begin{equation}
\label{e5.7}h_qh_p=h_p^2h_q^5,\qquad h_q^5h_p^2=h_ph_q\text{.} 
\end{equation}
The group $U$ acts on the set of these relations. Therefore the set $%
A=\{h_q,h_q^5\}$ is invariant with respect to $U$. By \ref{le24.3}) of Lemma 
\ref{le2.4}, $A\cdot A=\{h_q^2,h_q^4,h_q^6=id\}$ is also invariant. Let $%
\Gamma ^{\prime }$ be the subgroup of $\Gamma $ generated by $T_p$ and $%
\{id,h_q^2,h_q^4\}$. Then $\Gamma ^{\prime }$ is $U$-invariant, hence $%
G^{\prime }=\Gamma ^{\prime }U$ is a subgroup of $G$. Further, $[G:G^{\prime
}]=[\Gamma :\Gamma ^{\prime }]<\infty $. According to Proposition \ref{pp5.1}
the main group-theoretic lemma holds for $G^{\prime }$ (in this case $m=n=3$%
). We get a subgroup $U^{\prime }$ of $U$ such that $[U:U^{\prime }]<\infty $
and $U^{\prime }$ is a normal subgroup of $G^{\prime }$. Now $%
\bigcap\limits_{g\in G/G^{\prime }}gU^{\prime }g^{-1}$ is the desired
subgroup of $U$.\TeXButton{End Proof}{\endproof}

\begin{proposition}
\label{pp5.3}The main group-theoretic lemma holds for $n=2$ and arbitrary~$m$%
.
\end{proposition}

\TeXButton{Proof}{\proof}Let $\sigma \in T_p$, $\sigma \ne 1$. The group $U$
preserves $\sigma $. Therefore $\sigma T_q\sigma $ is invariant with respect
to $U$. Let $\Gamma ^{\prime }$ be the subgroup of $\Gamma $ generated by $%
T_q$ and $\sigma T_q\sigma $. Then $G^{\prime }=\Gamma ^{\prime }U$ is a
subgroup of $G$. We have $[G:G^{\prime }]=[\Gamma :\Gamma ^{\prime }]<\infty 
$. According to Proposition \ref{pp5.1} the main group-theoretic lemma holds
for $G^{\prime }$, so one can obtain the desired subgroup of $U$ just as in
the proof of Proposition \ref{pp5.2}.\TeXButton{End Proof}{\endproof}

\end{document}